\documentclass[11pt, reqno]{amsart}

\usepackage{amsfonts,amsmath,latexsym,amssymb,verbatim,amsbsy,amsthm,mathabx}
\usepackage{dsfont,bm}
\numberwithin{equation}{section}
\usepackage{enumitem}

\usepackage[margin=1in]{geometry}

\usepackage[dvipsnames]{xcolor}

\usepackage[colorlinks=true, pdfstartview=FitV, linkcolor=RoyalBlue,citecolor=ForestGreen, urlcolor=blue]{hyperref}

\definecolor{labelkey}{rgb}{0,0,1}

\theoremstyle{plain}
\newtheorem{THEOREM}{Theorem}[section]

\newtheorem{theorem}[THEOREM]{Theorem}

\newtheorem{lemma}[THEOREM]{Lemma}
\newtheorem{proposition}[THEOREM]{Proposition}

\theoremstyle{definition}

\newtheorem{definition}[THEOREM]{Definition}

\theoremstyle{remark}

\newcommand{\thm}[1]{Theorem~\ref{#1}}
\newcommand{\lem}[1]{Lemma~\ref{#1}}

\newcommand{\sect}[1]{Section~\ref{#1}}

\def \d {\delta}

\def \f {\varphi}

\def \n {\nabla}
\def \s {\sigma}

\def \D {\Delta}

\def \ba {{\bf a}}
\def \bb {{\bf b}}
\def \bc {{\bf c}}

\def \be {{\bf e}}

\def \bu {{\bf u}}
\def \bv {{\bf v}}
\def \bw {{\bf w}}
\def \bx {{\bf x}}
\def \by {{\bf y}}
\def \bz {{\bf z}}

\def \barPhi { \widebar{\Phi}}
\def \barx {\bar{x}}
\def \barv {\bar{v}}

\def \barX {\widebar{X}}
\def \barV {\widebar{V}}

\def \cD {\mathcal{D}}
\def \cE {\mathcal{E}}

\def \cG {\mathcal{G}}

\def \cK {\mathcal{K}}
\def \cL {\mathcal{L}}

\def \cP {\mathcal{P}}

\def \cR {\mathcal{R}}

\def \cV {\mathcal{V}}
\def \cW {\mathcal{W}}

\newcommand{\N}{\ensuremath{\mathbb{N}}}   
\newcommand{\R}{\ensuremath{\mathbb{R}}}   
\newcommand{\T}{\ensuremath{\mathbb{T}}}   

\newcommand{\E}{\ensuremath{\mathbb{E}}}

\def \one {{\mathds{1}}}

\renewcommand{\geq}{\geqslant}
\renewcommand{\ge}{\geqslant}
\renewcommand{\leq}{\leqslant}
\renewcommand{\le}{\leqslant}

\DeclareMathOperator{\supp}{supp} %
\DeclareMathOperator{\diam}{diam} %
\def \dx  {\, \mbox{d}x}

\def \dt  {\, \mbox{d}t}

\def \dy  {\, \mbox{d}y}

\def \dr  {\, \mbox{d}r}
\def \ds  {\, \mbox{d}s}

\def \dmu  {\, \mbox{d}\mu}

\def \dv  {\, \mbox{d}\bv}

\def \ddt  {\dfrac{\mbox{d\,\,}}{\mbox{d}t}}

\def \dd  {\mathrm{d}}

\def \ta {\tilde{a}}
\def \tb {\tilde{b}}

\def \bxd {\dot{\bx}}
\def \bvd {\dot{\bv}}

\def \lan {\langle}
\def \ran {\rangle}

\def \p {\partial}
\def \ra {\rightarrow}

\def \HI {H\"older inequality}

\renewcommand{\geq}{\geqslant}
\renewcommand{\ge}{\geqslant}
\renewcommand{\leq}{\leqslant}
\renewcommand{\le}{\leqslant}

\def \Lip {\mathrm{Lip}}

\def \P {\operatorname{P}}
\def \Var {\operatorname{Var}}

\def \dx {\dd\bx}
\def \dy {\dd \by}

\def \dyw {\dd\by \dd\bw}
\def \dXV {\dd X_0 \dd V_0}
\def \dnu {\dd \nu}
\def \dom {\T^d\times\R^{d}}

\begin{document}

\title[Mean-Field Limits of Deterministic and Stochastic Flocking Models]{Mean-Field Limits of Deterministic and Stochastic Flocking Models with Nonlinear Velocity Alignment}

\author{Vinh Nguyen$^\ast$}
\address{$^\ast$Michigan State University, Department of Mathematics,
East Lansing, MI 48824, USA}
\email{nguy1685@msu.edu }

\author{Roman Shvydkoy$^\dagger$}
\address{$^\dagger$University of Illinois at Chicago, Department of Mathematics, Statistics and Computer Science, Chicago,  IL 60607, USA}
\email{shvydkoy@uic.edu}

\author{Changhui Tan$^\ddagger$}
\address{$^\ddagger$University of South Carolina, Department of Mathematics, Columbia, SC 29208, USA}
\email{tan@math.sc.edu}


\date{\today}


\thanks{\textbf{Acknowledgment.}  
	 The work of R. Shvydkoy was supported in part by NSF grant DMS-2405326 and the Simons Foundation. 
	 The work of C. Tan was supported in part by NSF grant DMS-2238219.}

\begin{abstract}
We study the mean-field limit for a class of agent-based models describing flocking with nonlinear velocity alignment. Each agent interacts through a communication protocol $\phi$ and a non-linear coupling of velocities given by the power law $A(\bv) = |\bv|^{p-2}\bv$, $p > 2$. The mean-field limit is proved in two settings -- deterministic and stochastic. We then provide quantitative estimates on propagation of chaos for deterministic case in the case of the classical fat-tailed kernels, showing an improved convergence rate of the $k$-particle marginals to a solution of the corresponding Vlasov equation. The stochastic version is addressed with multiplicative noise depending on the local interaction intensity, which leads to the associated Fokker-Planck-Alignment equation. 

Our results extend the classical Cucker-Smale theory to the nonlinear framework which has received considerable attention in the literature recently.
\end{abstract}
\maketitle

\section{Introduction}\label{sec:intro}
The mathematical modeling of collective behavior, such as flocking of birds, schooling of fish, or swarming of bacteria, has been a fertile ground for interaction between mathematical analysis, probability theory, and statistical physics, see these comprehensive surveys of the subject \cite{albi2019vehicular, ben2005opinion, vicsek2012collective, motsch2014heterophilious,mucha2018cucker,shvydkoy2021dynamics}. A cornerstone in this field is the Cucker-Smale model \cite{cucker2007emergent}, a second-order particle system where agents align their velocities based on a weighted average of their neighbors' relative velocities. This model and its numerous variants have been shown to exhibit remarkable emergent properties, most notably flocking: the convergence of agents to a common velocity while maintaining a bounded spatial profile.

A fundamental question in the study of such interacting particle systems is their behavior as the number of particles $N$ tends to infinity. The formal limit is described by a kinetic equation, in which the probability distribution is transported along a force field generated by the distribution itself. The rigorous justification of this limit under the mean-field scaling and the associated propagation of chaos, which is the phenomenon in which particles become asymptotically independent as $N \to \infty$, are central problems in the mathematical physics and  modeling, see \cite{golse2016dynamics, jabin2014review, sznitman2006topics}.

There is a large and growing literature on mean-field limits and propagation of chaos for flocking models of Cucker-Smale type. 
For the classical Cucker-Smale system, a formal derivation of the associated kinetic equation via the BBGKY hierarchy was presented by Ha and Tadmor~\cite{ha2008particle}, while a rigorous justification through the mean-field limit was later established by Ha and Liu~\cite{ha2009simple}. The quantitative propagation of chaos for systems with fat-tailed communication protocols was subsequently investigated by Nguyen and Shvydkoy~\cite{nguyen2022propagation}, where explicit convergence rates in the Wasserstein-2 distance for finite marginals were obtained. 
Other works, including Natalini and Paul \cite{natalini2021mean} and related references therein, have also treated mean-field limits for the linear model under various kernel assumptions and interaction structures. 

While most classical results concern linear alignment, that is, interactions of the form $\bv_j - \bv_i$ weighted by a communication protocol, recent interest has shifted toward models with more singular or nonlinear interactions. These present significant mathematical challenges, but are often more biologically realistic or exhibit richer dynamics. 
\medskip

In this paper, we consider the following deterministic agent-based flocking model:
\begin{equation}\label{NVA}
    \begin{cases}
    \begin{aligned}
     \bxd_i &= \bv_i, \quad \bx_i(0) = \bx_i^0 \in \R^d\\
        \bvd_i &= \dfrac{1}{N}\sum\limits_{j= 1}^N \phi(\bx_i - \bx_j)A(\bv_j- \bv_i), \quad \bv_i(0) = \bv_i^0 \in \R^d.   
    \end{aligned}   
    \end{cases}
\end{equation}
Here $\phi$ is the \emph{communication protocol}, which measures the strength of the alignment interaction. We assume that $\phi$ is radially symmetric and non-increasing in the radial variable. A typical family of communication protocols we consider takes the form
\begin{equation}\label{phi}
 \phi(r) = (1+r)^{-\alpha},\quad \alpha\in[0,1), 	
\end{equation}
where $\phi$ is bounded, Lipschitz, and has a \emph{fat tail}, that is, there exists $r_0>0$ such that
\begin{equation*}
 \int_{r_0}^{\infty} \phi(r)\dr = \infty.
\end{equation*}

The mapping $A: \R^d \ra \R^d$ encodes the velocity coupling. For the linear choice $A(\bv) = \bv$, \eqref{NVA} reduces to the Cucker-Smale system. In this work we focus on the nonlinear mapping 
\begin{equation}\label{map-A}
    A(\bv) = |\bv|^{p-2}\bv, \quad p\ge 2.
\end{equation}
Such nonlinear velocity couplings were first introduced by Ha, Ha, and Kim \cite{ha2010emergent}. Subsequent works \cite{wen2012flocking,kim2020complete,black2024asymptotic} demonstrate that the choice of the nonlinear mapping has a decisive influence on the asymptotic flocking and alignment behavior, notably on the resulting convergence rates. In particular, Black and Tan \cite{black2024asymptotic} obtained quantitative bounds in which the rates depend explicitly on the parameters $p$ and $\alpha$. See Theorem \ref{thm:agent-flocking} for a detailed description.

The mean-field limit of \eqref{NVA} as $N\to\infty$ is given by the following Vlasov-type kinetic equation (see Definition \ref{def:ms-soln}):
\begin{equation}\label{VE}
    \begin{cases}
    \begin{aligned}
     &\p_t \mu + \bv\cdot \n_{\bx} \mu + \n_{\bv} \cdot (\mu F(\mu)) = 0, \quad \mu(0) = \mu_0 \in \cP(\R^{2d}),\\
     &F(\mu)(\bx,\bv,t) = \int_{\R^{2d}} \phi(\bx-\by) A(\bw-\bv) \dmu(\by,\bw). 
    \end{aligned}   
    \end{cases}
\end{equation}

The well-posedness theory for \eqref{VE} was established by Carrillo, Choi and Hauray \cite{carrillo2014local}. By formal hydrodynamic limits, one may further derive the macroscopic \emph{$p$-alignment system}:
\begin{equation}\label{pEA}
    \begin{cases}
    \begin{aligned}
     &\p_t \rho +  \n_{\bx} \cdot (\rho \bu) = 0, \quad \rho(0) = \rho_0: \R^{d} \to \R_+,\\
     &\p_t (\rho \bu) +\n_{\bx}\cdot (\rho \bu \otimes \bu) =  \int_{\R^{d}} \phi(\bx-\by) A\big(\bu(\by)-\bu(\bx)\big) \rho(\bx)\rho(\by) \dy. 
    \end{aligned}   
    \end{cases}
\end{equation}
We refer to the work of Tadmor \cite{tadmor2023swarming} and references therein for formal derivations and analysis of the $p$-alignment system. Remarkably, the systems \eqref{NVA}, \eqref{VE}, and \eqref{pEA} exhibit the same asymptotic flocking and alignment behavior (see also \cite{black2024asymptotic}).

The well-posedness theory of \eqref{pEA} and the rigorous derivation of the hydrodynamic limit were recently established by Black and Tan \cite{black2025hydrodynamic}. We also highlight the work of Choi, Fabisiak, and Peszek \cite{choi2025alignment}, who treated singular communication protocols of the form $\phi(r)=r^{-\alpha}$ with $\alpha\ge d$, establishing well-posedness and micro-to-macro mean-field limits.
\medskip

Our primary interest is the mean-field limit from the agent-based system \eqref{NVA} to the kinetic equation \eqref{VE}. Under the assumptions of smooth $\phi$ and the nonlinear mapping \eqref{map-A}, we derive a stability estimate analogous to the classical Dobrushin estimate \cite{dobrushin1979vlasov}: perturbations in the initial data, measured in the Wasserstein-1 distance, grow at most exponentially in time; see \eqref{W1:exp-bound}.

The kinetic equation \eqref{VE} inherits the same flocking and alignment estimates as the agent-based model \eqref{NVA}. Leveraging these estimates, we obtain enhanced stability bounds in which the growth is sub-exponential in time. As a consequence, we establish a rigorous mean-field limit for \eqref{NVA} and obtain enhanced quantitative propagation-of-chaos estimates in the Wasserstein-2 metric. Theorems \ref{thm:stability}--\ref{thm:chaos} present the detailed results. Our conclusions extend and generalize the quantitative estimates obtained in \cite{nguyen2022propagation} for the Cucker-Smale system with linear velocity coupling ($p=2$).

The kinetic equation \eqref{VE} enjoys the same flocking and alignment estimates as the agent-based model \eqref{NVA}. Using these estimates, we establish enhanced stability estimates where the growth is sub-exponential in time. Consequently, we establish a rigorous mean-field limit for agent-based system \eqref{NVA}, and derive enhanced quantitative propagation of chaos control, under Wasserstein-$2$ metric. See Theorems \ref{thm:stability}--\ref{thm:chaos} for detail descriptions, and Table \ref{table:rates} for explicit rates. 
Our results broaden the quantitative framework developed in \cite{nguyen2022propagation} for the Cucker-Smale system with linear velocity coupling ($p=2$), extending it to the nonlinear regime.

In the stochastic setting, a number of works addressed the classical linear alignment force with multiplicative or additive noise. The first result in this direction belongs to  Bolley, Ca\~nizo, and Carrillo \cite{bolley2011stochastic}, in the case of constant noise and linear alignment force. Strength-dependent noise was treated in \cite{shvydkoy2024environmental}.  Choi and Salem \cite{choi2019cucker} analyzed stochastic particle systems with multiplicative noise depending on velocity, deriving stochastic mean-field limits and phase transition phenomena. Friesen and Kutoviy \cite{friesen2020stochastic} considered jump-type stochastic interactions and proved the propagation of chaos through McKean-Vlasov formulations. Those studies treat additive or velocity-dependent noise.

We consider the following stochastic agent-based flocking model:
\[
\dd \bv_i = \dfrac{1}{N}\sum\limits_{j= 1}^N \phi(\bx_i - \bx_j)A(\bv_j- \bv_i) \dt + \sqrt{2h(s_{i})}\dd W_i,
\]
incorporating multiplicative noise whose amplitude depends on the local interaction strength $s_i$. We develop an analogous mean-field theory in this stochastic setting and prove convergence of the empirical measure to the associated Fokker-Planck-Alignment equation \eqref{FPE}.

Overall, our results significantly broaden the scope of existing work, which has largely been confined to linear velocity couplings or simpler noise structures. They provide a unified quantitative framework for both deterministic and stochastic flocking models with nonlinear velocity alignment.
\medskip

The remainder of the paper is organized as follows. 
\sect{sec:prelim} introduces notation, recalls the flocking estimates for the agent-based system, presents key auxiliary lemmas, and states the main results. 
In \sect{sec:MFL}, we study the deterministic flocking model with nonlinear alignment and establish global well-posedness together with the mean-field limit. 
\sect{sec:chaos} is devoted to quantitative propagation-of-chaos estimates, where explicit Wasserstein-$2$ convergence rates are derived for systems with fat-tailed communication protocols. 
In \sect{sec:StochasticMFL}, we extend the analysis to the stochastic model with locally dependent multiplicative noise and prove convergence of the empirical measure to the corresponding Fokker-Planck-Alignment equation. 
Finally, \sect{sec:appendix} contains the proofs of the key auxiliary lemmas used throughout the analysis.
\subsection*{Notations}
We use $C$ to denote a positive constant whose value may change from line to line. 
When two distinct constants appear in the same expression, we write $C$ and $\bar{C}$. 
We often employ the Japanese bracket $\lan t\ran := \sqrt{1+t^{2}}$. 
The functions $a(t)$ and $b(t)$ will serve as auxiliary quantities. 
The notation $A \lesssim B$ means that there exists a positive constant $C$ such that $A \leq C B$.
\section{Preliminaries and Statement of Main Results}\label{sec:prelim}
In this section we collect preliminary material and state the main results of the paper. We first introduce notation and recall basic properties of the Wasserstein distance, which will serve as the principal metric throughout the analysis. We then review the flocking and alignment estimates for the agent-based model with nonlinear velocity alignment under fat-tailed communication protocols. Next, we present a set of auxiliary lemmas for systems of differential inequalities, which provide the quantitative bounds needed in the subsequent stability and propagation-of-chaos arguments. The section concludes with the statements of the main results, including the Dobrushin-type stability estimate for the kinetic equation, the deterministic mean-field limit, quantitative propagation-of-chaos estimates, and the stochastic mean-field limit.

\subsection{Wasserstein distance} We first recall the definition and key properties of the Wasserstein distance (see \cite{ambrosio2021lectures,golse2016dynamics}), a central metric for measuring the proximity of probability measures in our analysis.

Throughout, denote by $\cP(\R^{2d})$ the space of Borel probability measures on $\R^{2d}$, and for $m>0$ set
\begin{equation*}
    \cP_{m}(\R^{2d}):= \left\{\mu \in \cP(\R^{2d}) : \int_{\R^{2d}} |\bv|^{m}\dmu(\bx,\bv) <\infty\right\}.
\end{equation*}

Given $\mu,\nu \in \cP_m(\R^d)$, we write $\Pi(\mu,\nu)$ for the set of couplings between $\mu$ and $\nu$; in other words, $\Pi(\mu,\nu)$ consists of all Borel probability measures $\pi$ on $\R^d \times \R^d$ whose first and second marginals are $\mu$ and $\nu$, respectively. Equivalently, a measure $\pi \in \cP(\R^d \times \R^d)$ belongs to $\Pi(\mu,\nu)$ if and only if  
\[
\iint_{\R^d \times \R^d} \bigl(\varphi(\bx) + \psi(\by)\bigr)\, \dd\pi(\bx,\by) 
= \int_{\R^d} \varphi(\bx)\, \dmu(\bx) + \int_{\R^d} \psi(\by)\, \dnu(\by)
\]  
for every pair of functions $\varphi,\psi \in C(\R^d)$ with at most polynomial growth, i.e., $\varphi(\bz),\psi(\bz) = O(|\bz|^m)$ as $|\bz|\to\infty$.

\begin{definition}
For $m \geq 1$ and $\mu,\nu \in \cP_m(\R^d)$, the Wasserstein-$m$ distance (also known as the Monge--Kantorovich distance) between $\mu$ and $\nu$, denoted by $\cW_m(\mu,\nu)$, is defined by  
\[
\cW_m(\mu,\nu) := \inf_{\pi \in \Pi(\mu,\nu)} 
\Big( \iint_{\R^d \times \R^d} |\bx-\by|^m \, \dd\pi(\bx,\by) \Big)^{1/m}.
\]
\end{definition}
For a function $\varphi: \R^d \to \R$, the Lipschitz constant of $\varphi$, denoted by $\Lip(\varphi)$, is defined by
\[
\Lip(\varphi):=\sup_{\bx\ne \by\in\R^d}\frac{|\varphi(\bx)-\varphi(\by)|}{|\bx-\by|}.
\]
The set of Lipschitz functions on $\R^d$, denoted by $\Lip(\R^d)$, consists all functions whose Lipschitz constant is finite. 

The Wasserstein-$1$ distance has the following dual formulation:
\begin{equation}\label{W1dual}
\cW_1(\mu,\nu)=\sup_{\substack{\varphi\in\Lip(\R^d)\\ \Lip(\varphi)\le 1}}\left|\int_{\R^d}\varphi(\bz)\dmu(\bz)-\int_{\R^d}\varphi(\bz)\dnu(\bz)\right|.
\end{equation}

\subsection{Flocking and alignment with fat-tailed communication protocol}
Let us recall the flocking and alignment properties of the agent-based model \eqref{NVA}.  
For $t\ge0$, define the spatial and velocity diameters
\begin{equation*}
    \cD(t) =  \max_{i,j=1,\ldots,N} |\bx_i(t) - \bx_j(t)|,\quad \cV(t) = \max_{i,j=1,\ldots,N} |\bv_i(t) - \bv_j(t)|.
\end{equation*}
The flocking and alignment behaviors can be interpreted by
\[
\sup_{t\geq0} \cD(t)<\infty,\quad\text{and}\quad \lim_{t\to\infty}\cV(t)=0,
\]
respectively.

For the nonlinear alignment model \eqref{NVA}, the evolution of $(\cD,\cV)$ is known to satisfy the following closed system of differential inequalities \cite{ha2010emergent}:
\begin{equation}\label{SDDI}
    \begin{cases}
      \cD'(t) \leq \cV(t),\\[2mm]
      \cV'(t)\leq - c_p \phi(\cD(t)) \cV(t)^{p-1}, \quad c_p = 2^{2-p}.
    \end{cases}
\end{equation}
Further analysis of the paired inequalities show flocking and alignment phenomena, provided that the communication protocol $\phi$ is fat-tailed, satisfying \eqref{fat-tail}.

In \cite{black2024asymptotic}, quantitative bounds were obtained for the family of communication protocols in \eqref{phi}, or more generally, $\phi$ is bounded, Lipschitz, and there exist $\lambda, \Lambda > 0$, $\alpha \in [0,1)$ such that
    \begin{equation}\label{fat-tail}
        \lambda r^{-\alpha} \leq \phi(r) \leq \Lambda r^{-\alpha}, \quad \forall r > r_0.
    \end{equation}
Sharp asymptotic bounds for $\cD(t)$ and $\cV(t)$ were derived. These rates depend explicitly on the nonlinear exponent $p$ and the tail parameter $\alpha$. We summarize the result below.    
\begin{theorem}[Flocking and alignment estimates \cite{black2024asymptotic}]\label{thm:agent-flocking}
Suppose $(\cD(t),\cV(t))$ satisfy the paired inequality \eqref{SDDI} with bounded initial data $(\cD_0,\cV_0)$, and the communication protocol $\phi$ satisfies \eqref{fat-tail}.
Then, for any $t\geq0$, we have the following bounds:
    \begin{enumerate}[label = \textnormal{(\roman*).}]
        \item if $ p \in (2,3)$, then there exists $D < \infty$ such that
        \begin{equation*}\label{flocking-fat-3-}
          \sup_{t\geq 0}\cD(t)  = D, \quad \cV(t) \lesssim \langle t\rangle^{-\frac{1}{p-2}}; 
        \end{equation*}
        \item if $p > 3$, then
        \begin{equation*}\label{flocking-fat-3+}
            \cD(t)\lesssim \lan t\ran^{1- \frac{1-\alpha}{p-\alpha-2}}, \quad \cV(t) \lesssim  \lan t\ran^{- \frac{1-\alpha}{p-\alpha-2}}.
        \end{equation*}
        \item if $p = 3$, then
        \begin{equation*}\label{flocking-fat-3}
            \cD(t)\lesssim \big(\log\lan t\ran\big)^{\frac{1}{1-\alpha}}, \quad \cV(t) \lesssim  \lan t\ran^{-1} \big(\log\lan t\ran\big)^{\frac{\alpha}{1-\alpha}};
        \end{equation*}
    \end{enumerate}   
\end{theorem}
In contrast to the linear velocity coupling case ($p=2$), where $\cV(t)$ decays exponentially fast, the nonlinear regime $p>2$ yields only polynomial decay. Moreover, flocking (i.e., bounded $\cD(t)$) holds for all $p<3$, while for $p\ge3$ the spatial diameter may grow in time, though only at a sublinear rate.
\subsection{Key auxiliary lemmas}
We shall require a collection of technical estimates for certain closed systems of differential inequalities sharing the same structural form as \eqref{SDDI}. These systems admit quantitative bounds analogous to those in Theorem~\ref{thm:agent-flocking}. The resulting estimates will be used repeatedly in the sequel, in particular for deriving quantitative stability and propagation-of-chaos bounds. Proofs of the auxiliary lemmas stated below are deferred to \sect{sec:appendix}.
%
%
\begin{lemma}\label{lem:ab-aux}
 Given $\beta \in [0,+\infty)$, suppose that $(a(t), b(t))$ is a nonnegative solution to the system:
 \begin{equation}\label{ab-aux-sys}
     \begin{cases}
         a'(t) \leq b(t), \quad a(0) = a_0\in \R_+,\\[2mm]
         b'(t) \leq C\langle t\rangle ^{-\beta}a(t) + g(t), \quad b(0) = b_0\in \R_+, \quad t\geq 0,
     \end{cases}
 \end{equation}
 where $C> 0$ is a constant, and $g(t)$ is a given nonnegative source term. Then the following holds:
 \begin{enumerate}[label = \textnormal{(\roman*).}]
     \item if $\beta > 2$, letting
     \[\cG_1(t):=\int_0^t g(s)\ds,\]
     then there exist positive constants $C_1$ and $C_2$ such that
     \begin{equation}\label{eq:aux1}
    a(t) \leq a_0 + C_1 t + C_2\int_0^t \cG_1(s)\ds,\quad b(t) \leq C_1 + C_2 \cG_1(t),  \quad t\geq 0.
     \end{equation}
     \item if $\beta \in [0,2)$, letting $\gamma := 1- \frac{\beta}{2}$, $\bar C:= \frac{1-\gamma + \sqrt{(1-\gamma)^2 + 4C}}{2\gamma}$, and
     \begin{equation*}
      \cG_2(t) := \int_0^t \lan s\ran^{1-\gamma} \,e^{-\bar C \,\lan s\ran^{\gamma}} g(s) \ds,
     \end{equation*}
      then there exist positive constants $C_3$ and $C_4$ such that
     \begin{equation}\label{eq:aux2}
     a(t) \leq \big(C_1 + C^{-1}\bar C \gamma\, \cG_2(t)\big) e^{\bar C \langle t\rangle^{\gamma}}, \quad b(t) \leq \big(C_2 + \cG_2(t)\big) \langle t\rangle^{-(1-\gamma)}\,e^{\bar C \langle t\rangle^{\gamma}}, \quad t\geq 0.
     \end{equation}
     \item If $\beta = 2$, letting $\zeta:= (1+\sqrt{1+4C})/2$ and
     \[
     \cG_3(t):= \int_0^t \lan s\ran^{-(\zeta-1)}g(s) \ds,
     \]
     then there exist positive constants $C_5$ and $C_6$ such that 
     \[
     a(t) \leq \big(C_5+C^{-1}\zeta\cG_3(t)\big) \langle t\rangle^{\zeta}, \quad b(t) \leq \big(C_6+\cG_3(t)\big) \langle t\rangle^{\zeta -1}, \quad t\geq 0.
     \]
 \end{enumerate}
\end{lemma}
We briefly comment on the behavior of \eqref{ab-aux-sys} in the absence of a source term, i.e., when $g(t)\equiv0$.  
For $\beta=0$, part (ii) recovers the classical exponential growth bound
\[
a(t)\leq C_1 e^{\bar Ct},\quad b(t)\leq C_2 e^{\bar Ct}.
\]
When $\beta\in(0,2)$, the solution exhibits sub-exponential but super-polynomial growth.  
The case $\beta=2$ is critical and leads to polynomial bounds.  
Finally, for $\beta>2$, the growth of $a(t)$ is at most linear, while $b(t)$ remains uniformly bounded.

In Theorem~\ref{thm:agent-flocking}(iii), corresponding to the critical case $p=3$, the flocking and alignment estimates involve a logarithmic correction. Accordingly, in the critical regime $\beta=2$, we require a variant of \eqref{ab-aux-sys} that incorporates an additional logarithmic factor. The corresponding estimate is given next.
\begin{lemma}\label{lem:ab-aux2}
Given $\alpha\in[0,1)$, suppose that $(a(t), b(t))$ is a nonnegative solution to the system:
 \begin{equation}\label{ab-aux-sys2}
     \begin{cases}
         a'(t)\leq b(t), \quad a(0) = a_0\in \R_+,\\[2mm]
         b'(t) \leq C\langle t\rangle ^{-2} (\log \lan t\ran )^{\frac{2\alpha}{1-\alpha}} a(t) + g(t), \quad b(0) = b_0\in \R_+, \quad t\geq 0,
     \end{cases}
 \end{equation}
 where $C> 0$ is a constant, and $g(t)$ is a given nonnegative source term. Letting $\theta :=\frac{1}{1-\alpha}$, $\bar C := \frac{1 +\sqrt{1 + 4C}}{2\theta}$, and
\[
  \cG_4(t) := \int_0^t \langle s\rangle(\log \langle s\rangle)^{-(\theta-1)}\,e^{-\bar C (\log \langle s\rangle)^{\theta}}g(s) \ds,
\]
  then there exist positive constants $C_7$ and $C_8$ such that
\begin{equation}\label{eq:aux3}
     a(t) \leq \big(C_7 + C^{-1}\bar C\theta \cG_4(t)\big)\, e^{\bar C(\log \langle t\rangle)^{\theta}}, \quad b(t) \leq \big(C_8 + \cG_4(t)\big)\, \langle t\rangle^{-1}(\log \langle t\rangle)^{\theta-1}\,e^{\bar C (\log \langle t\rangle)^{\theta}}, \quad t\geq 0.
\end{equation}
\end{lemma}
Lemma~\ref{lem:ab-aux2} provides the refined estimates required to treat the critical case, including the corresponding sharp rate. We note that the presence of the logarithmic factor $(\log\langle t\rangle)^{\frac{2\alpha}{1-\alpha}}$ in \eqref{ab-aux-sys2} leads to a bound that is no longer polynomial in time, in contrast to Lemma~\ref{lem:ab-aux}(iii), except in the special case $\alpha=0$.
\subsection{Statement of main results}
In this subsection, we state our main results together with several accompanying remarks. 

We first establish a Dobrushin-type stability estimate for the nonlinear kinetic equation \eqref{VE}, which forms the backbone of the subsequent mean-field and propagation-of-chaos analysis.

\begin{definition}\label{def:ms-soln}
A map $\mu: [0,T) \ra \cP_{2}(\R^{2d}), \; t\mapsto \mu_t$, is called a \emph{measure-valued solution} to \eqref{VE} with initial data $\mu_0$ if it satisfies the following conditions:
\begin{itemize}
    \item[(i).] $\mu$ is weakly* continuous, 
\item[(ii).] For any $\f \in C^{\infty}_0([0,T)\times \R^{2d})$ and $0<t<T$,
\begin{align*}
   \int_{\R^{2d}} \f(t,\bx,\bv)\dmu_t(\bx,\bv)  = &     \int_{\R^{2d}}\f(0,\bx,\bv) \dmu_0(\bx,\bv) \\
   & \quad\quad  + \int_0^t \int_{\R^{2d}} \big[\p_s\f + \bv\cdot\n_\bx\f +F(\mu_s)\cdot\n_\bv\f  \big]\dmu_s(\bx,\bv)\ds.   
\end{align*}
\end{itemize}
\end{definition}

\begin{theorem}[Stability]\label{thm:stability}
Let $\mu$ and $\nu$ be two measure-valued solutions of \eqref{VE} on a common time interval of existence $[0,T)$, corresponding to initial data $\mu_0$ and $\nu_0$, respectively. 
Assume that the initial measures have compact supports contained in a common compact subset $\Omega \subset \R^{2d}$, that is,
\[
\supp \mu_0 \cup \supp \nu_0 \subseteq \Omega.
\]
Then there exist constants $C,\bar{C}>0$, depending only on $\Omega$ and $\phi$, such that
\begin{equation}\label{W1:exp-bound}
    \cW_1(\mu_t,\nu_t) \le C e^{\bar{C} t}\,\cW_1(\mu_0,\nu_0),
    \qquad \forall\, t \in [0,T).
\end{equation}
Moreover, if the communication protocol $\phi$ satisfies the fat-tailed condition \eqref{fat-tail}, we have:
\begin{enumerate}[label=\textnormal{(\roman*).}]
    \item if $p \in (2,3)$, then
    \[
    \cW_1(\mu_t,\nu_t) \leq C\lan t\ran\log\lan t\ran\cW_1(\mu_0,\nu_0);
    \]
    \item if $p > 3$,  then 
    \[
     \cW_1(\mu_t,\nu_t) \leq C\lan t\ran^{\frac{(1+\alpha)(p-3)}{2(p-\alpha-2)}}e^{\bar C\lan t\ran^{\frac{(1+\alpha)(p-3)}{2(p-\alpha-2)}}}\cW_1(\mu_0,\nu_0);
    \]
    \item if $p=3$, then 
    \[
    \cW_1(\mu_t,\nu_t) \leq C(\log \langle t\rangle)^{\frac{1}{1-\alpha}}e^{\bar C (\log \langle t\rangle)^{\frac{1}{1-\alpha}}}\cW_1(\mu_0,\nu_0),
    \]
\end{enumerate}
where $C, \bar{C}$ are positive constants which depend only on initial support, $\phi$ and $p$.
\end{theorem}
Building on this stability estimate, we establish well-posedness of the kinetic equation \eqref{VE} and provide a rigorous justification of the mean-field limit from the agent-based model \eqref{NVA}.
\begin{theorem}[Mean-field limit]\label{thm:MFL}
Given any measure $\mu_0 \in \cP_{2}(\R^{2d})$ with compact support, there exists a unique measure-valued solution to \eqref{VE} with the initial condition $\mu_0$. More specifically, this solution is the weak limit of the empirical measure  built on the solution to the agent-based system \eqref{NVA}.      
\end{theorem}
\medskip

Our next main result is the quantitative propagation-of-chaos estimate under fat-tailed communication. We start with introducing the notations. 

Suppose the initial measure $\mu_0$ is absolutely continuous with respect to the Lebesgue measure,
\begin{equation*}
 \dmu_0(\bx,\bv) = f_0(\bx,\bv)\dx\dv,   
\end{equation*}
where $f_0:\R^{2d}\ra \R_+$ is a density function, then the solution of \eqref{VE} is of the form $\dmu_t = f(\bx,\bv,t)\dx\dv$ with density $f$ satisfying
\begin{equation}\label{VEfa}
    \begin{cases}
    \begin{aligned}
     & \p_t f + \bv\cdot \n_{\bx} f + \n_{\bv} \cdot (f F(f)) = 0, \quad f(0) = f_0,\\
     & F(f)(\bx,\bv,t) = \int_{\R^{2d}} \phi(\bx-\by) A(\bw-\bv) f(\by,\bw,t) \dyw.
    \end{aligned}   
    \end{cases}
\end{equation}
Denote by $f^N: \R^{2dN}\ra \R_+$ the $N$-particle density function, which solves the Liouville equation  
\begin{equation}\label{LEa}
    \begin{cases}
    \begin{aligned}
     & \p_t f^{N} + \sum_{i=1}^N \bv_i \cdot \n_{\bx_i} f^N + \sum_{i=1}^N \n_{\bv_i} \cdot(f^N F_i^N) = 0, \quad f^N(0) = f_0^{\otimes N},\\
     & F_i^N(\bx_1,\ldots,\bx_N,\bv_1,\ldots,\bv_N) = \frac1N \sum_{j=1}^N \phi(\bx_i - \bx_j)A(\bv_j - \bv_i).
    \end{aligned}   
    \end{cases}
\end{equation}
Due to the symmetries of the initial data and the forces, the solution will remain symmetric with respect to permutations of pairs $(\bx_i,\bv_i)$ for all time.
We further define the $k$-th marginal by
\begin{equation*}
f^{(k)}_t(\bx_1,\bv_1,\ldots,\bx_k,\bv_k) = \int_{\R^{2d(N-k)}} f^N(\bx_1,\bv_1,\ldots,\bx_N,\bv_N,t)\,\dd\bx_{k+1}\dd\bv_{k+1}\ldots \dd\bx_N\dd\bv_N.
\end{equation*}

\begin{theorem}[Propagation of chaos]\label{thm:chaos}
Suppose $\phi$ satisfies the fat-tailed condition \eqref{fat-tail} and
 $f_0\in C^1_0(\R^{2d})$ is a compactly supported initial distribution. Let $f, f^N$ be the solutions to \eqref{VEfa} and \eqref{LEa}, respectively. Then for all $N\in \N$, $k < N$, and $t\geq 0$,
\begin{enumerate}[label=\textnormal{(\roman*).}]
     \item if $p\in (2,3)$ then
     \begin{equation*}
 \cW_2(f^{(k)}_t,f^{\otimes k}_t) \leq C \sqrt{k}\min\left\{1,\frac{t}{\sqrt{N}}\right\};
\end{equation*}
\item if $p > 3$ then
     \[
      \cW_2(f^{(k)}_t,f^{\otimes k}_t) \leq C\sqrt{\frac{k}{N}}e^{\bar C \langle t\rangle^{\frac{(1+\alpha)(p-3)}{2(p-\alpha-2)}}};
     \]     
\item if $p=3$ then
     \[
\cW_2(f^{(k)}_t,f^{\otimes k}_t) \leq C\sqrt{\frac{k}{N}}(\log \langle t\rangle)^{\frac{\alpha}{1-\alpha}}e^{\bar C(\log \langle t\rangle)^{\frac{1}{1-\alpha}}},
     \]
 \end{enumerate}   
where constants $C,\bar{C}$ depend only on $p$, $\diam(\supp f_0)$ and $\phi$.
\end{theorem}

\renewcommand{\arraystretch}{1.5}
\begin{table}[h]
 \caption{Quantitative estimates on flocking models with nonlinear velocity alignment}\label{table:rates}
 \begin{tabular}{lcccc}
    \hline
	& $p=2$ & $2<p<3$ & $p=3$ & $p>3$\\ \hline
	Flocking:\,\,$\cD(t)$ & $D$ & $D$ & $\big(\log\lan t\ran\big)^{\frac{1}{1-\alpha}}$ &  $\lan t\ran^{1- \frac{1-\alpha}{p-\alpha-2}}$\\ \hline
	Alignment:\,\,$\cV(t)$ & $e^{-\kappa t}$ & $\langle t\rangle^{-\frac{1}{p-2}}$ & $\lan t\ran^{-1} \big(\log\lan t\ran\big)^{\frac{\alpha}{1-\alpha}}$ & $\lan t\ran^{- \frac{1-\alpha}{p-\alpha-2}}$\\ \hline 	
	Stability:\,\,$\cW_1(\mu_t,\nu_t)$ \hspace{-.1in} & \multicolumn{2}{c}{$\lan t\ran\log\lan t\ran$} & $(\log \langle t\rangle)^{\frac{1}{1-\alpha}}e^{\bar C (\log \langle t\rangle)^{\frac{1}{1-\alpha}}}$ & $\lan t\ran^{\frac{(1+\alpha)(p-3)}{2(p-\alpha-2)}}e^{\bar C\lan t\ran^{\frac{(1+\alpha)(p-3)}{2(p-\alpha-2)}}}$\\[5pt] \hline 
	PoC:\,\,$\cW_2(f^{(k)}_t,f^{\otimes k}_t)$ & \multicolumn{2}{c}{$\sqrt{k}\min\left\{1,\frac{t}{\sqrt{N}}\right\}$} & $\sqrt{\frac{k}{N}}(\log \langle t\rangle)^{\frac{\alpha}{1-\alpha}}e^{\bar C(\log \langle t\rangle)^{\frac{1}{1-\alpha}}}$& $\sqrt{\frac{k}{N}}e^{\bar C \langle t\rangle^{\frac{(1+\alpha)(p-3)}{2(p-\alpha-2)}}}$ \\[5pt] \hline 
 \end{tabular}
\end{table}
Table~\ref{table:rates} summarizes the quantitative bounds on stability and the corresponding propagation-of-chaos (PoC) estimates. We remark that in the regime $p\in(2,3)$, the flocking property is sufficient to recover the same rate for the PoC bound as in the classical Cucker-Smale model with linear velocity coupling ($p=2$), despite the lack of exponential alignment. For $p>3$, the resulting bounds grow super-polynomially but remain sub-exponential in time. In the borderline case $p=3$, the estimates contain an additional logarithmic correction.
\medskip

Our final result focus on the stochastic model with strength-dependent noise
\begin{equation}\label{SNVA}
    \begin{cases}
    \begin{aligned}
     \dd \bx_i &= \bv_i\dt,\quad (\bx_i,\bv_i) \in \T^d\!\times\!\R^{d},\\
        \dd \bv_i &= \dfrac{1}{N}\sum\limits_{j= 1}^N \phi(\bx_i - \bx_j)A(\bv_j- \bv_i) \dt + \sqrt{2h(s_{i})}\dd W_i,   
    \end{aligned}   
    \end{cases}
\end{equation}
where $W_i$'s are independent Brownian motions in $\R^d$.
The function $s_i$ is called the \emph{strength function}, defined by
\[
s_i := \dfrac{1}{N}\sum\limits_{j= 1}^N \phi(\bx_i - \bx_j),
\]
and the mapping $h:\R\to \R$ is Lipschitz and copositive, which means $h(r) > 0$ if $r>0$. Denote by $\T^d$ the $d-$dimensional torus. We consider the following stochastic system:

We will show that the stochastic mean field limit for this system is
\begin{equation}\label{FPE}
\p_t f + \bv\cdot \n_{\bx} f = h(s_\rho)\Delta_v f - \n_{\bv} \cdot (f F(f)) ,
\end{equation}
where $F(f)$ given by 
\[
F(f)(\bx,\bv) = \int_{\dom} \phi(\bx-\by) A(\bw-\bv) f(\by,\bw,t) \dd\by \dd\bw.
\]

For stochastic setting, the following result will be proved in \sect{sec:StochasticMFL}.
\begin{theorem}\label{thm:StochasticMFL}
 Suppose that $\{(\bx_i, \bv_i)\}_{i=1}^N$ is a solution to \eqref{SNVA} with joint distribution law $f^N$ such that $f^N_0= f_0^{\otimes N}$. Let $f$ solve \eqref{FPE} on $[0,T]$ with initial data $f_0$. If $f$ satisfies
 \begin{align}
   &\sup_{t\in [0,T]}\phi\ast\rho(\bx) \ge \underline{\rho}>0,\label{pos}\\
  & \sup_{t\in [0,T]}  \int_{\dom} e^{c_0|\bv|^{p-1}} f(\bx,\bv,t)\dd\bx\dd\bv < \infty  \text{ for some } c_0 > 0. \label{moment}
 \end{align}
then we have the following mean-field approximation of $f$ in law: for any $\f \in \Lip$,
 \begin{equation*}
    \E \Big|  \frac1N \sum_{i=1}^N \varphi(\bx_i(t),\bv_i(t)) - \int_{\dom} \varphi(\bx,\bv) f(t,\bx,\bv) \dd\bx\dd\bv \Big|^2  \leq \frac{C}{N^{e^{-Ct}}} \quad \forall t \in [0,T],
 \end{equation*}
and propagation of chaos
 \begin{equation*}
 W_2^2(f^{(k)}, f^{\otimes k}) \leq \frac{Ck}{N^{e^{-Ct}}} \quad \forall\, t \in [0,T].
  \end{equation*}
\end{theorem}

\section{ Mean-Field Limit For Deterministic Model}\label{sec:MFL}
This section is dedicated to establishing the mean-field limit for the deterministic flocking model with nonlinear velocity alignment \eqref{NVA}. We begin by formally associating the Vlasov-type kinetic equation \eqref{VE} with \eqref{NVA} via characteristic paths. We then analyze the kinetic equation \eqref{VE}, showing that it inherits the same flocking and alignment behavior as the agent-based model \eqref{NVA}.

The core of the section is the derivation of a Dobrushin-type stability estimate in the Wasserstein-$1$ distance, demonstrating that the flow generated by \eqref{VE} is Lipschitz continuous with respect to its initial data. We then combine this estimate with the flocking and alignment control to obtain an enhanced stability bound. This stability result serves as the key ingredient in the rigorous mean-field analysis, allowing us to show that the empirical measure of the agent-based system converges to the solution of the kinetic equation as the number of agents tends to infinity.

\subsection{Flocking and alignment estimates}
Suppose that $\mu: [0,T) \ra \cP_{2}(\R^{2d})$ is a measure-valued solution to \eqref{VE} with initial data $\mu_0$. It is known from the optimal transport theory that $\mu$ is a push-forward of $\mu_0$ along the characteristics $(X_\mu, V_\mu)$:  
\begin{equation}\label{eq:XV}
\begin{cases}
\begin{aligned}
    &\ddt X_\mu(t,\bx,\bv) = V_\mu(t,\bx,\bv),\quad X_\mu(0,\bx,\bv) = \bx,\\
    & \ddt V_\mu(t,\bx,\bv) = \int_{\R^{2d}} \phi(X_\mu' - X_\mu)A(V_\mu' - V_\mu) \dmu_0(\bx',\bv'), \quad V_\mu(0,\bx,\bv) =\bv.
\end{aligned}
\end{cases}
\end{equation}
Here and in the following, we write $X'_\mu = X_\mu(t,\bx',\bv')$ and $V'_\mu = V_\mu(t,\bx',\bv')$.  
For simplicity, we also denote $\bz = (\bx,\bv)$ and $\bz' = (\bx',\bv')$. If we assume that the support of the initial measure $\mu_0$ is contained in a compact set $\Omega$, then we define
\begin{equation}\label{def:DV-kin}
    \cD_\Omega(t) = \max_{\bz', \bz'' \in \Omega} |X_\mu(t,\bz') - X_\mu(t,\bz'')|, \quad \cV_{\Omega}(t) = \max_{\bz',\bz'' \in \Omega} |V_\mu(t,\bz') - V_\mu(t,\bz'')|.
\end{equation}

The following proposition shows that $(\cD_\Omega,\cV_\Omega)$ satisfies the differential inequalities \eqref{SDDI}. The proof is analogous to that in \cite{ha2010emergent,black2024asymptotic}.

\begin{proposition}
	Suppose $(X_\mu, V_{\mu})$ satisfies \eqref{eq:XV}. Then, $(\cD_\Omega, \cV_\Omega)$ defined in \eqref{def:DV-kin} satisfies the differential inequalities \eqref{SDDI}, that is, for almost all $t\ge0$, 
	\begin{equation}\label{kin-flock-sys}
    \begin{cases}
        \cD_{\Omega}'(t) \leq \cV_{\Omega}(t)\\[2mm]
        \cV_\Omega'(t)\leq - c_p \phi(\cD_\Omega(t)) \cV(t)^{p-1}, \quad c_p = 2^{2-p}.
    \end{cases}
\end{equation}
\end{proposition}
\begin{proof}
 From \eqref{eq:XV}$_1$ and Rademacher's lemma, we immediately get
 \[
 \ddt \cD_{\Omega}(t) \leq \cV_{\Omega}(t).
 \]
 For the second inequality in \eqref{kin-flock-sys}, we fix a time $t$ and take $\bz', \bz'' \in \Omega$ to be the maximizing points, such that
 \[
\cV_{\Omega}(t) = |V_\mu(t,\bz') - V_\mu(t,\bz'')|.
\]
We remark that the selections $\bz'$ and $\bz''$ depend on time, and are in general neither unique nor continuous with respect to time.
From \eqref{eq:XV}$_2$ and Rademacher's lemma we obtain
\begin{align*}
   \frac12 \ddt \cV_\Omega^2(t) =&\, (V'_\mu-V_\mu'')\cdot \ddt (V'_\mu-V_\mu'')\notag\\
   \leq & \,\phi(\cD_\Omega(t)) \int_{\R^{2d}} (V'_\mu-V_\mu'')\cdot\Big[|V_\mu-V_\mu'|^{p-2}(V_\mu-V_\mu') - |V_\mu-V_\mu''|^{p-2} (V_\mu-V_\mu'') \Big]\dmu_0(\bz)\notag\\
   \leq & - c_p \phi(\cD_\Omega(t)) \int_{\R^{2d}}|V'_\mu-V_\mu''|^{p} \dmu_0(\bz) = - c_p \phi(\cD_\Omega(t)) \cV_\Omega(t)^{p},
\end{align*}
where in the penultimate step, we have used an elementary inequality
\[
 (V'_\mu-V_\mu'')\cdot\Big[|V_\mu-V_\mu'|^{p-2}(V_\mu-V_\mu') - |V_\mu-V_\mu''|^{p-2} (V_\mu-V_\mu'') \Big] \leq -2^{2-p} |V'_\mu-V_\mu''|^{p},
\]
which can be justified by applying Lemma \ref{lem:eleineq} below with $\ba = V_\mu'$, $\bb = V_\mu''$, and $\bc = V_\mu$. 

Collecting the estimates above, we conclude with the desired inequalities \eqref{kin-flock-sys}.
\end{proof}

\begin{lemma}\label{lem:eleineq}
 Let $p\geq2$. Then for any $\ba, \bb, \bc \in\R^d$,
 \begin{equation}\label{eq:eleineq}
  (\ba - \bb)\cdot\big( |\ba-\bc|^{p-2}(\ba-\bc) - |\bb-\bc|^{p-2}(\bb-\bc) \big)\geq 2^{2-p}|\ba-\bb|^p,
 \end{equation}
 where the equality is attained when $\ba=\bb$ or $\bc=\frac{\ba+\bb}{2}$.
\end{lemma}

\begin{proof}
When $\ba=\bb$, it is obvious that \eqref{eq:eleineq} holds with an equality.  Now, assume $\ba\neq\bb$.
By rotation, translation, and scaling invariance, the inequality \eqref{eq:eleineq} holds for $(\ba, \bb, \bc)$ if and only if it holds for $(\widetilde\ba, \widetilde\bb, \widetilde\bc)$ with
\[
\widetilde\ba = \tfrac{2}{|\ba-\bb|}\,\cR_{\ba-\bb}\left(\ba - \tfrac{\ba+\bb}{2}\right) = \be_1,\quad
\widetilde\bb = \tfrac{2}{|\ba-\bb|}\,\cR_{\ba-\bb}\left(\bb - \tfrac{\ba+\bb}{2}\right) = -\be_1,\quad
\widetilde\bc = \tfrac{2}{|\ba-\bb|}\,\cR_{\ba-\bb}\left(\bc - \tfrac{\ba+\bb}{2}\right),
\]
where $\cR_{\ba-\bb}$ denotes the rotation transformation such that $\cR_{\ba-\bb}(\ba-\bb) = |\ba-\bb|\,\be_1$.
Therefore, it suffices to show that for any $\widetilde\bc\in\R^d$,
 \[
 2\be_1\cdot \big( |\be_1-\widetilde\bc|^{p-2}(\be_1-\widetilde\bc) - |-\be_1-\widetilde\bc|^{p-2}(-\be_1-\widetilde\bc) \big)\geq 2^{2-p}\cdot 2^p = 4.
 \] 
 Express $\widetilde\bc = (\widetilde c_1, \widetilde\bc_r)$, where $\widetilde c_1$ is the first component of $\widetilde\bc$, and $\widetilde\bc_r\in\R^{d-1}$ consists the remaining components. Then, we have
 \begin{align*}
  & \be_1\cdot \big( |\be_1-\widetilde\bc|^{p-2}(\be_1-\widetilde\bc) - |-\be_1-\widetilde\bc|^{p-2}(-\be_1-\widetilde\bc) \big)\\
  & = \big((1-\widetilde c_1)^2+|\widetilde\bc_r|^2\big)^{\frac{p-2}{2}}(1-\widetilde c_1) + \big((1+\widetilde c_1)^2+|\widetilde\bc_r|^2\big)^{\frac{p-2}{2}}(1+\widetilde c_1)\\
  & \geq |1-\widetilde c_1|^{p-2}(1-\widetilde c_1) + |1+\widetilde c_1|^{p-2}(1+\widetilde c_1)\geq2,
 \end{align*}
 where the two inequalities attain equality when $\widetilde\bc_r=\mathbf{0}$ and $\widetilde c_1=0$, respectively.
 We conclude with the inequality \eqref{eq:eleineq}, where the equality is attained when $\widetilde\bc =\mathbf{0}$, or equivalently $\bc = \frac{\ba+\bb}{2}$.
\end{proof}

We remark that taking $\bx = \ba - \bc$ and $\by = \bb-\bc$, we deduce from \eqref{eq:eleineq} the following interesting inequality:
\begin{equation}\label{eq:eleineq2}
(\bx - \by) \cdot (|\bx|^{p-2}\bx - |\by|^{p-2}\by) \ge 2^{2-p} |\bx-\by|^p \ge 0,\quad \forall\,\bx,\by\in\R^d.
\end{equation}

From \eqref{kin-flock-sys}$_2$, we have $\cV_\Omega'(t)\le 0$, which yields a maximum principle for the velocity diameter:
\begin{equation}\label{gen-p-velocity}
    \cV_{\Omega}(t) \le \cV_{\Omega}(0)
    \quad \forall\, t\ge 0.
\end{equation}
Substituting this bound into \eqref{kin-flock-sys}$_1$ leads to a linear growth estimate for the spatial diameter:
\begin{equation*}
    \cD_{\Omega}(t) \le C(1+t)
    \quad \forall\, t\ge 0,
\end{equation*}
where $C$ is a constant depending only on the $\cD_\Omega(0)$ and $\cV_\Omega(0)$.

Since the system \eqref{kin-flock-sys} has the same structural form as \eqref{SDDI}, the enhanced flocking and alignment estimates follow directly from Theorem~\ref{thm:agent-flocking} under the assumption that the communication protocol is fat-tailed. We summarize the resulting bounds below.
\begin{theorem}\label{thm:kin-flocking}
Assume that the communication protocol $\phi$ satisfies the fat-tailed condition \eqref{fat-tail}. 
Suppose that $\mu \in C_{w^*}(\R^+; \cP_{2}(\R^{2d}))$ is a measure-valued solution to \eqref{VE} with compactly supported initial data $\mu_0$, and let $\Omega \subset \R^{2d}$ be a compact set with $\supp \mu_0 \subseteq \Omega$. 
Then the following holds for all $t\geq 0$:
\begin{enumerate}[label=\textnormal{(\roman*).}]
        \item if $ p \in (2,3)$ then  
        \[
          \sup_{t>0}\cD_\Omega(t)  = D <\infty, \quad \cV_\Omega(t) \lesssim \langle t\rangle^{-\frac{1}{p-2}};
        \]
        \item if $p > 3$, then
        \[
            \cD_\Omega(t)\lesssim \lan t\ran^{1- \frac{1-\alpha}{p-\alpha-2}}, \quad \cV_\Omega(t) \lesssim  \lan t\ran^{- \frac{1-\alpha}{p-\alpha-2}}.
        \]
        \item if $p = 3$, then
        \[
            \cD_{\Omega}(t)\lesssim \big(\log\lan t\ran\big)^{\frac{1}{1-\alpha}}, \quad \cV_\Omega(t) \lesssim  \lan t\ran^{-1} \big(\log\lan t\ran\big)^{\frac{\alpha}{1-\alpha}};
        \]
    \end{enumerate} 
\end{theorem}

\subsection{Stability}
Based upon the flocking and alignment estimates, we now prove our first main result, \thm{thm:stability}.

Denote by $(X_\mu,V_\mu)$ and $(X_\nu,V_\nu)$ the characteristic flows associated with $\mu$ and $\nu$, respectively.  
Then, by the dual formulation of the Wasserstein-$1$ distance \eqref{W1dual}, we have
\begin{align}\label{W1}
\cW_1(\mu_t,\nu_t)
&= \sup_{\substack{\varphi\in \Lip(\R^{2d})\\ \Lip(\varphi)\le 1}}
\left|\int_{\R^{2d}}\varphi(\bz)\,\dmu_t(\bz)
- \int_{\R^{2d}}\varphi(\bz)\,\dnu_t(\bz)\right|\notag\\
&= \sup_{\substack{\varphi\in \Lip(\R^{2d})\\ \Lip(\varphi)\le 1}}
\left|\int_{\R^{2d}}\varphi(X_\mu,V_\mu)\,\dmu_0(\bz)
- \int_{\R^{2d}}\varphi(X_\nu,V_\nu)\,\dnu_0(\bz)\right|\notag\\
&\le \sup_{\substack{\varphi\in \Lip(\R^{2d})\\ \Lip(\varphi)\le 1}}
\left|\int_{\R^{2d}}\varphi(X_\mu,V_\mu)\,\dmu_0(\bz)
- \int_{\R^{2d}}\varphi(X_\mu,V_\mu)\,\dnu_0(\bz)\right|\notag\\
&\quad + \sup_{\substack{\varphi\in \Lip(\R^{2d})\\ \Lip(\varphi)\le 1}}
\left|\int_{\R^{2d}}\big[\varphi(X_\mu,V_\mu)-\varphi(X_\nu,V_\nu)\big]\dnu_0(\bz)\right| \notag\\
&\le \big(\|\nabla X_\mu\|_\infty + \|\nabla V_\mu\|_\infty\big)
\,\cW_1(\mu_0,\nu_0)
+ \int_{\R^{2d}}|(X_\mu,V_\mu)-(X_\nu,V_\nu)|\,\dnu_0(\bz)\notag\\
&\le \big(\|\nabla X_\mu\|_\infty + \|\nabla V_\mu\|_\infty\big)\cW_1(\mu_0,\nu_0)
+ \|X_\mu-X_\nu\|_\infty + \|V_\mu-V_\nu\|_\infty.
\end{align}
The next step is to estimate the quantities 
$\|\nabla X_\mu\|_\infty + \|\nabla V_\mu\|_\infty$ and 
$\|X_\mu-X_\nu\|_\infty + \|V_\mu-V_\nu\|_\infty$.
From equation~\eqref{eq:XV}$_1$, we have
\[
\ddt\nabla X_\mu = \nabla V_\mu.
\]
By evaluating this equality at points where 
$\|\nabla X_\mu\|_\infty$ and $\|X_\mu-X_\nu\|_\infty$ are achieved, 
and invoking Rademacher's lemma, we obtain
\begin{equation}\label{X}
\ddt\|\nabla X_\mu\|_\infty \le \|\nabla V_\mu\|_\infty, 
\qquad
\ddt\|X_\mu-X_\nu\|_\infty \le \|V_\mu-V_\nu\|_\infty.
\end{equation}
Next, by differentiating \eqref{eq:XV}$_2$ we find
\begin{align*}
\ddt\nabla V_\mu 
&= \int_{\R^{2d}} 
\nabla\phi(X'_\mu - X_\mu)\,(\nabla^T X_\mu)\otimes A(V'_\mu-V_\mu)\dmu_0(\bz') \nonumber\\
& \hspace{-.2in} - \int_{\R^{2d}} \phi(X'_\mu - X_\mu)\Big[ |V'_\mu-V_\mu|^{p-2}\mathbf{I} 
+ (p-2)|V'_\mu-V_\mu|^{p-4}(V'_\mu-V_\mu)\otimes(V'_\mu-V_\mu)\Big]
\nabla V_\mu\dmu_0(\bz'),
\end{align*}
where $\mathbf{I}$ denotes the $d$-by-$d$ identity matrix.
Since the matrix inside the bracket of the second integral is positive definite, 
we deduce that
\begin{equation}\label{der-gradV}
\ddt\|\nabla V_\mu\|_\infty 
\le \|\nabla\phi\|_\infty\,\|\nabla X_\mu\|_\infty\,\cV_{\Omega}^{p-1}(t).
\end{equation}
Then, for a general communication protocol $\phi$, using \eqref{gen-p-velocity} we obtain 
\begin{equation*}
\ddt\|\nabla V_\mu\|_\infty \le C\|\nabla X_\mu\|_\infty,
\end{equation*}
 where $C$ is a positive constant depending on $\phi$ and initial data.
Combining this with the first inequality in~\eqref{X} and applying Gr\"onwall's inequality yields, for all $t\in[0,T)$,
\begin{equation}\label{gen-phi-grad}
\|\nabla X_\mu\|_\infty + \|\nabla V_\mu\|_\infty \le C e^{Ct}.
\end{equation}
If $\phi$ satisfies fat-tailed condition \eqref{fat-tail}, then combining \eqref{der-gradV} and the estimates in Theorem \ref{thm:kin-flocking} we obtain for all $t\in [0,T)$:
\begin{enumerate}[label=\textnormal{(\roman*).}]
    \item if $p \in (2,3)$ then
    \begin{equation*}\label{gradV3-}
\ddt\|\nabla V_\mu\|_\infty 
\le C\|\nabla X_\mu\|_\infty\,\langle t\rangle^{-\frac{p-1}{p-2}}.
\end{equation*}
\item  if $p>3$ then
    \begin{equation*}\label{gradV3+}
\ddt\|\nabla V_\mu\|_\infty 
\le C\|\nabla X_\mu\|_\infty\,\lan t\ran^{- \frac{(1-\alpha)(p-1)}{p-\alpha-2}}.
\end{equation*}
\item  if $p =3$ then
    \begin{equation*}\label{gradV3}
\ddt\|\nabla V_\mu\|_\infty 
\le C\|\nabla X_\mu\|_\infty\,\lan t\ran^{-2} \big(\log\lan t\ran\big)^{\frac{2\alpha}{1-\alpha}}.
\end{equation*}
\end{enumerate}
Here $C>0$ denoting a general constant depending only on $\phi$, $p$ and the initial data, it is different for each case of $p$.  
Set
\[
a(t):=\|\nabla X_\mu\|_\infty,\qquad 
b(t):=\|\nabla V_\mu\|_\infty.
\]
Then, the dynamics of $(a,b)$ satisfies the paired inequalities in \eqref{ab-aux-sys} with $g(t)\equiv0$, and the parameter
\begin{equation}\label{eq:beta}
\beta := \begin{cases}
 	\,\,\frac{p-1}{p-2}>2 & \text{for}\,\, p\in(2,3),\\[2mm]
 	\,\,\frac{(1-\alpha)(p-1)}{p-\alpha-2} \in(0,2)\quad & \text{for}\,\, p>3.
 \end{cases}
\end{equation}
For the borderline case $p=3$, the dynamics satisfies \eqref{ab-aux-sys2} with $g(t)\equiv0$.
Applying the key Lemmas \ref{lem:ab-aux} and \ref{lem:ab-aux2}, we arrive at the following estimates.
\begin{lemma}\label{claim:grad-est}
For all $t\in [0,T)$, we have
\begin{enumerate}[label=\textnormal{(\roman*).}]
    \item if $p\in (2,3)$ then
    \[
	\|\nabla X_\mu\|_\infty \lesssim \lan t\ran,\quad
    \|\nabla V_\mu\|_\infty\lesssim 1; 
    \]
   \item if $p>3$ then 
    \[
         \|\nabla X_\mu\|_\infty \lesssim e^{\bar C\lan t\ran^{\gamma}}, \quad  
\|\nabla V_\mu\|_\infty \lesssim \lan t\ran^{-(1-\gamma)} e^{\bar C\lan t\ran^{\gamma}},
    \]
 where $\gamma = 1- \frac\beta2 = \frac{(1+\alpha)(p-3)}{2(p-\alpha-2)}$;
    \item if $p=3$ then
    \[
\|\nabla X_\mu\|_\infty  \lesssim \, e^{\bar C(\log \langle t\rangle)^{\frac{1}{1-\alpha}}}, \quad \|\nabla V_\mu\|_\infty \lesssim \, \langle t\rangle^{-1}(\log \langle t\rangle)^{\frac{\alpha}{1-\alpha}}\,e^{\bar C (\log \langle t\rangle)^{\frac{1}{1-\alpha}}}.
    \] 
\end{enumerate}  
\end{lemma}
We now estimate the time derivative of $\|V_\mu - V_\nu\|_{\infty}$. 
Evaluating at a maximizing point $\bz$, we obtain
\begin{align}
\dfrac{1}{2}\ddt  \|V_\mu - V_\nu\|^2_{\infty} 
&\leq (V_\mu - V_\nu)\cdot \Bigg(
\int_{\R^{2d}} \phi(X'_\mu - X_\mu)A(V'_\mu - V_\mu) \,\dmu_0(\bz') \notag\\
&\hspace{4.5cm}-\int_{\R^{2d}} \phi(X'_\nu - X_\nu)A(V'_\nu - V_\nu) \,\dnu_0(\bz') \Bigg)\notag\\[1ex]
&\leq (V_\mu - V_\nu)\cdot \int_{\R^{2d}} \phi(X'_\mu - X_\mu)A(V'_\mu - V_\mu)\big( \dmu_0(\bz') - \dnu_0(\bz')\big)\notag\\
&\quad + (V_\mu - V_\nu)\cdot \int_{\R^{2d}} \big(\phi(X'_\mu - X_\mu)-\phi(X'_\nu - X_\nu)\big)A(V'_\mu - V_\mu) \,\dnu_0(\bz')\notag\\
&\quad + \int_{\R^{2d}} \phi(X'_\nu - X_\nu)(V_\mu - V_\nu)\cdot \big(A(V'_\mu - V_\mu) - A(V'_\nu - V_\nu)\big)\,\dnu_0(\bz').\notag
\end{align}
The last term is negative due to Lemma \ref{lem:eleineq}. Therefore,
\begin{align}
\ddt \|V_\mu - V_\nu\|_{\infty}
&\le C\|\phi\|_{W^{1,\infty}(\Omega)}
\Big(\|\nabla X_\mu\|_{\infty}\cV_{\Omega}^{p-1} + \|\nabla V_\mu\|_{\infty}\cV_{\Omega}^{p-2}\Big)\cW_1(\mu_0,\nu_0)\notag\\
&\quad + C\|\nabla \phi\|_{\infty}\|X_\mu - X_\nu\|_{\infty}\cV_{\Omega}^{p-1}\label{eq:der-distV}
\end{align}
Hence, for general $\phi$ we use the bounds from \eqref{gen-p-velocity} and \eqref{gen-phi-grad} to have
\begin{equation*}
 \ddt \|V_\mu - V_\nu\|_{\infty} \leq Ce^{Ct}\cW_1(\mu_0,\nu_0) + C\|X_\mu - X_\nu\|_{\infty}.   
\end{equation*}
Combining this inequality with the second estimate in \eqref{X}, we arrive at
\[
\ddt\big(\|X_\mu - X_\nu\|_{\infty} + \|V_\mu - V_\nu\|_{\infty}\big)
\le C\big(\|X_\mu - X_\nu\|_{\infty} + \|V_\mu - V_\nu\|_{\infty}\big)
+ Ce^{Ct}\cW_1(\mu_0,\nu_0).
\]
Applying Gr\"onwall's inequality yields
\begin{equation}\label{eq:exp-stab}
\|X_\mu - X_\nu\|_{\infty} + \|V_\mu - V_\nu\|_{\infty}
\le Ce^{Ct}\cW_1(\mu_0,\nu_0), \qquad \forall\, t\in [0,T),
\end{equation}
where $C>0$ depend on the initial data and the interaction potential $\phi$.  
Substituting \eqref{eq:exp-stab} and \eqref{gen-phi-grad} into \eqref{W1}, we obtain the exponential growth stability estimate for a general kernel $\phi$:
\begin{equation}\label{W1-gen-phi}
\cW_1(\mu_t,\nu_t) \le Ce^{Ct}\cW_1(\mu_0,\nu_0), \qquad \forall\, t\in [0,T),  
\end{equation}
which is the conclusion \eqref{W1:exp-bound}.

\smallskip
In the case $\phi$ satisfies the fat-tailed condition \eqref{fat-tail}, utilizing the  estimates in Theorem \ref{thm:kin-flocking} and Lemma \ref{claim:grad-est}, from \eqref{eq:der-distV} we obtain:
\begin{enumerate}[label=\textnormal{(\roman*).}]
    \item if $p\in (2,3)$ then
\begin{align*}
\ddt \|V_\mu - V_\nu\|_{\infty} &
\le C\big(\lan t\ran \cdot \lan t\ran^{-\frac{p-1}{p-2}} + 1\cdot\lan t\ran^{-1}\big)\,\cW_1(\mu_0,\nu_0) + C \|X_\mu - X_\nu\|_{\infty}\,\lan t\ran^{-\frac{p-1}{p-2}}\\
& \leq C\lan t\ran^{-\beta}\|X_\mu - X_\nu\|_{\infty} + C\lan t\ran^{-1}\cW_1(\mu_0,\nu_0),
\end{align*}
where we have used the fact $1-\frac{p-1}{p-2}<-1$ and the definition of $\beta$ in \eqref{eq:beta};
\item if $p>3$, then
\begin{align*}
\ddt \|V_\mu - V_\nu\|_{\infty} &
\le C\big(e^{\bar C\lan t\ran^{\gamma}}\cdot\lan t\ran^{-\frac{(1-\alpha)(p-1)}{p-\alpha-2}} + \lan t\ran^{-(1-\gamma)}e^{\bar C\lan t\ran^{\gamma}}\cdot\lan t\ran^{-\frac{(1-\alpha)(p-2)}{p-\alpha-2}}\big)\,\cW_1(\mu_0,\nu_0)\\
&\quad + C \|X_\mu - X_\nu\|_{\infty}\,\lan t\ran^{-\frac{(1-\alpha)(p-1)}{p-\alpha-2}}\\
& \leq C\lan t\ran^{-\beta}\|X_\mu - X_\nu\|_{\infty} + C\lan t\ran^{-\beta}e^{\bar C\lan t\ran^{\gamma}}\cW_1(\mu_0,\nu_0);
\end{align*}
where we have used the fact $-\frac{(1-\alpha)(p-1)}{p-\alpha-2}>-(1-\gamma)-\frac{(1-\alpha)(p-2)}{p-\alpha-2}$ and the definition of $\beta$;
\item if $p= 3$, then
\begin{align*}
\ddt \|V_\mu - V_\nu\|_{\infty} &
\le C \|X_\mu - X_\nu\|_{\infty}\,\lan t\ran^{-2} \big(\log\lan t\ran\big)^{\frac{2\alpha}{1-\alpha}} + C\Big(e^{\bar C(\log \langle t\rangle)^{\frac{1}{1-\alpha}}}\cdot\lan t\ran^{-2} \big(\log\lan t\ran\big)^{\frac{2\alpha}{1-\alpha}}\\
& \qquad\qquad+ \langle t\rangle^{-1}(\log \langle t\rangle)^{\frac{\alpha}{1-\alpha}}\,e^{\bar C (\log \langle t\rangle)^{\frac{1}{1-\alpha}}}\cdot\lan t\ran^{-1} \big(\log\lan t\ran\big)^{\frac{\alpha}{1-\alpha}}\Big)\,\cW_1(\mu_0,\nu_0)\\
& \leq C\lan t\ran^{-2} \big(\log\lan t\ran\big)^{\frac{2\alpha}{1-\alpha}}\|X_\mu - X_\nu\|_{\infty} +C \langle t\rangle^{-2}(\log \langle t\rangle)^{\frac{2\alpha}{1-\alpha}}\,e^{\bar C (\log \langle t\rangle)^{\frac{1}{1-\alpha}}}\cW_1(\mu_0,\nu_0).
\end{align*}
\end{enumerate}
Set 
\[a(t):=\|X_\mu - X_\nu\|_{\infty},\qquad b(t):=\|V_\mu - V_\nu\|_{\infty}.\] 
Then, the dynamics of $(a, b)$ satisfies the paired inequalities \eqref{ab-aux-sys} (or \eqref{ab-aux-sys2} for the borderline case) with  $a_0 = 0, b_0 = 0$, and with source term $g(t)$. Applying the key Lemmas \ref{lem:ab-aux} and \ref{lem:ab-aux2}, we obtain the following bounds.
\begin{lemma}\label{claim:dist-est}
For all $t\in [0,T)$,
\begin{enumerate}[label=\textnormal{(\roman*).}]
    \item if $p\in (2,3)$ then 
    \[
    \|X_\mu - X_\nu\|_{\infty}\lesssim t\log\lan t\ran\cW_1(\mu_0,\nu_0),\quad
    \|V_\mu - V_\nu\|_{\infty} \lesssim \log\lan t\ran\cW_1(\mu_0,\nu_0);
    \]
    \item if $p > 3$ then
    \[
     \|X_\mu - X_\nu\|_{\infty}\lesssim \lan t\ran^{\gamma}e^{\bar C\lan t\ran^{\gamma}}\cW_1(\mu_0,\nu_0),\quad 
     \|V_\mu - V_\nu\|_{\infty}\lesssim \lan t\ran^{2\gamma-1}e^{\bar C\lan t\ran^{\gamma}}\cW_1(\mu_0,\nu_0),
    \]
    where $\gamma = \frac{(1+\alpha)(p-3)}{2(p-\alpha-2)}$;
   \item if $p = 3$ then
   \[
   \|X_\mu - X_\nu\|_{\infty}\lesssim (\log \langle t\rangle)^{\theta}e^{\bar C (\log \langle t\rangle)^{\theta}}\cW_1(\mu_0,\nu_0),\quad
	\|V_\mu - V_\nu\|_{\infty} \lesssim \langle t\rangle^{-1}(\log \langle t\rangle)^{2\theta-1} e^{\bar C (\log \langle t\rangle)^{\theta}} \cW_1(\mu_0,\nu_0),
	\] 
	where $\theta = \frac{1}{1-\alpha}$.
\end{enumerate}    
\end{lemma}

\begin{proof}
For (i), we have $g(t) = C\lan t\ran^{-1}\cW_1(\mu_0,\nu_0)$. Compute
\[
      \cG_1(t) = C\cW_1(\mu_0,\nu_0)\int_0^t \lan s\ran^{-1} \ds = C\log\lan t\ran \cW_1(\mu_0,\nu_0).
\]
Then we deduce from \eqref{eq:aux1} that
\[
a(t)\lesssim \int_0^t\cG_1(s)\ds\lesssim t\log\lan t\ran\cW_1(\mu_0,\nu_0),\quad
b(t)\lesssim \cG_1(t)\lesssim \log\lan t\ran\cW_1(\mu_0,\nu_0).
\]
For (ii), we have $g(t) = C\lan t\ran^{-\beta}e^{\bar C\lan t\ran^{\gamma}}\cW_1(\mu_0,\nu_0).$ Compute
\[
      \cG_2(t) = C\cW_1(\mu_0,\nu_0)\int_0^t \lan s\ran^{\gamma-1} \ds \leq C\gamma^{-1} \cW_1(\mu_0,\nu_0) \lan t\ran^{\gamma}.
\]
Then, we deduce from \eqref{eq:aux2} that
\[
 a(t)\lesssim \cG_2(t) e^{\bar C \lan t\ran^\gamma} \lesssim \lan t\ran^{\gamma}e^{\bar C\lan t\ran^{\gamma}}\cW_1(\mu_0,\nu_0),\quad 
 b(t)\lesssim \cG_2(t) \lan t\ran^{-(1-\gamma)} e^{\bar C \lan t\ran^\gamma} \lesssim \lan t\ran^{2\gamma-1}e^{\bar C\lan t\ran^{\gamma}}\cW_1(\mu_0,\nu_0).
\]
For (iii), we have $g(t) = C \langle t\rangle^{-2}(\log \langle t\rangle)^{\frac{2\alpha}{1-\alpha}}\,e^{\bar C (\log \langle t\rangle)^{\frac{1}{1-\alpha}}}\cW_1(\mu_0,\nu_0)$. Compute
\[
	\cG_4(t) = C \cW_1(\mu_0,\nu_0)\int_0^t \lan s\ran^{-1}(\log \lan s\ran)^{\theta-1}\ds = C\theta^{-1} (\log \lan t\ran)^{\theta}\cW_1(\mu_0,\nu_0).
\]
Then, we deduce from \eqref{eq:aux3} that
\begin{align*}
 a(t) & \lesssim \cG_4(t)\, e^{\bar C(\log \langle t\rangle)^{\theta}} \lesssim (\log \lan t\ran)^{\theta}e^{\bar C(\log \langle t\rangle)^{\theta}}\cW_1(\mu_0,\nu_0),\\
 b(t) & \lesssim \cG_4(t)\, \langle t\rangle^{-1}(\log \langle t\rangle)^{\theta-1}\,e^{\bar C (\log \langle t\rangle)^{\theta}} \lesssim \lan t\ran^{-1}(\log \langle t\rangle)^{2\theta-1}\,e^{\bar C (\log \langle t\rangle)^{\theta}}\cW_1(\mu_0,\nu_0).
\end{align*}
\end{proof}

\medskip
Substituting the bounds from Lemmas \ref{claim:grad-est} and \ref{claim:dist-est} into \eqref{W1}, we obtain the following stability estimates for the fat-tailed $\phi$:
\begin{enumerate}[label=\textnormal{(\roman*).}]
    \item if $p \in (2,3)$ then
    \[
    \cW_1(\mu_t,\nu_t) \lesssim \lan t\ran\log\lan t\ran\cW_1(\mu_0,\nu_0);
    \]
    \item if $p > 3$ then 
    \[
     \cW_1(\mu_t,\nu_t) \lesssim\lan t\ran^{\gamma}e^{\bar C\lan t\ran^{\gamma}}\cW_1(\mu_0,\nu_0),\quad \gamma=\frac{(1+\alpha)(p-3)}{2(p-\alpha-2)};
    \]
    \item if $p=3$, then there exist constants $C,\bar{C}$ satisfying
    \[
    \cW_1(\mu_t,\nu_t) \lesssim (\log \langle t\rangle)^{\theta}e^{\bar C (\log \langle t\rangle)^{\theta}}\cW_1(\mu_0,\nu_0), \quad \theta = \frac{1}{1-\alpha}.
    \]
\end{enumerate}
Together with \eqref{W1-gen-phi}, we conclude the proof of Theorem \ref{thm:stability}.

\subsection{Mean-field Limit}
With the stability of the mean-field dynamics established, we now prove the existence and uniqueness of solutions to the kinetic equation \eqref{VE} by taking the limit of the particle system. The strategy is to consider a sequence of empirical measures $\mu_t^N$ corresponding to the particle system \eqref{NVA} and apply the Dobrushin stability estimate. This estimate implies that the sequence is Cauchy in the Wasserstein-$1$ distance, ensuring its convergence to a limit $\mu_t$. We then verify that this limit satisfies the weak formulation of the Vlasov equation, thus constituting the desired mean-field solution.

\begin{proof}[Proof of \thm{thm:MFL}]
For any $N\in \N$, choose $(\bx_k^0, \bv^0_k) \in \supp \mu_0, \, k = 1, \ldots, N$ such that
\[
\mu_0^N : = \frac{1}{N}\sum_{k=1}^N \d_{\bx_k^0}\otimes\d_{\bv_k^0}  \overset{*}{\rightharpoonup}\mu_0\quad  \text{ as } N \ra \infty.
\]

Define the empirical measures
\[
\mu_t^N := \frac{1}{N}\sum_{k= 1}^N  \d_{\bx_k(t)}\otimes\d_{\bv_k(t)},
\]
where $(\bx_k(t),\bv_k(t))$ is the solution to \eqref{NVA}.
Testing with $\varphi \in C_0^\infty([0,T)\times \R^{2d})$, we have that $\mu^N$ is a measure-valued solution to \eqref{VE} with initial data $\mu_0^N$. Thus, applying Theorem \ref{thm:stability} there exist constants $C,c>0$ such that
\[
\cW_1(\mu_t^{N}, \mu_t^{M}) \leq Ce^{cT} \cW_1(\mu_0^{N}, \mu_0^{M}), \quad \text{for } N, M>0, \ t\leq T.
\]
Hence $\{\mu_t^{N}\}_N$  is weakly$^*$-Cauchy in the Banach space $(\cP_2(\R^{2d}), \cW_1)$, and hence converges to a limit $\mu_t\in \cP_2(\R^{2d})$. Moreover,
\begin{equation*}
\cW_1(\mu_t^{N}, \mu_t) \leq C_T \cW_1(\mu_0^{N}, \mu_0), \quad \text{for } N>0, \ t\leq T.
\end{equation*}
Next we will prove that the map $t \to \mu_t$ is weak$^*$-continuous. Firstly, we note that for $\psi \in C^\infty_0 (\R^{2d})$ the sequence $\big\{\int_{\R^{2d}} \varphi(\bx,\bv)\dmu^N_{t}(\bx,\bv)\big\}_N$ is uniformly Lipschitz continuous on $[0, T]$. Indeed, for $t\in [0, T)$ and $\D t>0$ with $t+\D t \in [0,T]$ we have
\begin{align*}
    &\left|\int_{\R^{2d}} \psi(\bx,\bv)\dmu^N_{t+\D t}(\bx,\bv)-\int_{\R^{2d}} \psi(\bx,\bv)\dmu^N_{t}(\bx,\bv) \right | \\
    &\hspace{1.2in} \le   \int_{\R^{2d}} \left|\psi (X_{\mu^N} (t+\D t),V_{\mu^N} (t+\D t)) -\psi(X_{\mu^N}(t),V_{\mu^N}(t))\right| \dmu^N_0 (\bx,\bv)\\
     &\hspace{1.2in} \le |\n \psi|_{\infty}\int_{\R^{2d}}\Big(\big|X_{\mu^N} (t+\D t) - X_{\mu^N} (t)\big|+  \big|V_{\mu^N} (t+\D t) - V_{\mu^N}(t)\big| \Big)  \dmu^N_0 (\bx,\bv)\\
     & \hspace{1.2in} \le C \D t.
\end{align*}
For the last inequality we used the uniform Lipschitzness of $\{X_{\mu^N}\}_N, \{V_{\mu^N}\}_N$ on $[0,T]$. Then, letting $N\ra +\infty$ we obtain
\[
\left|\int_{\R^{2d}} \psi(\bx,\bv)\dmu_{t+\D t}(\bx,\bv)-\int_{\R^{2d}} \psi(\bx,\bv)\dmu_{t}(\bx,\bv) \right | \le  C \D t, 
\]
which implies the weak$^*$-continuity of the map $t \to \mu_t$.

\smallskip
We will show that this $\mu$ is a measure-valued solution to \eqref{VE} with the given initial $\mu_0$. Because $\mu^N$ is a measure-valued solution, for any test function $\f \in C^{\infty}_0([0,T)\times\R^{2d})$,
\begin{equation*}
\begin{split}
 \int_{\R^{2d}} \f(t,\bx,\bv)\dmu^N_t(\bx,\bv) =&     \int_{\R^{2d}}\f(0,\bx,\bv) \dmu^N_0(\bx,\bv) \\
 & \qquad + \int_0^t \int_{\R^{2d}} \big[\p_s\f + \bv\cdot\n_\bx\f +F(\mu^N_s)\cdot\n_\bv\f  \big]\dmu^N_s(\bx,\bv)\ds.   
\end{split}
\end{equation*}
All linear terms weakly converge to the natural limits. For the nonlinear term, we note that for all $t\in [0,T)$ and $N\in \N$, there exists $R_T>0$ such that $\supp \mu_t^N \subseteq B_{R_T}(0)$, a ball in $\R^{2d}$ with radius $R_T$. The family of functions $\big\{\phi(\bx-\cdot) A(\cdot-\bv)\big\}_{(\bx,\bv)}$ indexing by $(\bx,\bv)\in B_{R_T}(0)$ is uniformly Lipschitz on $B_{R_T}(0)$ with the common Lipschitz constant denoted by $L_{R_T}$. Thus, it is precompact in $C(B_{R_T}(0))$. Therefore, $F(\mu_t^N)(\bx,\bv) \ra F(\mu_t(\bx,\bv)$ uniformly on $B_{R_T}(0)$. We have
\[
\begin{split}
 &\Big| \int_0^t \int_{\R^{2d}} F(\mu^N_s)(\bx,\bv)\cdot\n_\bv\f(\bx,\bv)  \dmu^N_s(\bx,\bv)\ds - \int_0^t \int_{\R^{2d}} F(\mu_s)(\bx,\bv)\cdot\n_\bv\f(\bx,\bv)  \dmu_s(\bx,\bv)\ds \Big| \\
 &\quad\leq \Big| \int_0^t \int_{\R^{2d}} F(\mu^N_s)(\bx,\bv)\cdot\n_\bv\f(\bx,\bv)  \big(\dmu^N_s(\bx,\bv)-\dmu_s(\bx,\bv)\big)\ds\Big|\\
 &\quad \quad + \Big| \int_0^t \int_{\R^{2d}} \big(F(\mu^N_s)(\bx,\bv) - F(\mu_s)(\bx,\bv)\big)\cdot\n_\bv\f(\bx,\bv)  \dmu_s(\bx,\bv)\ds  \Big|\\
 &\quad\leq \int_0^t L_{R_T} |\nabla\varphi|_{\infty}\cW_1(\mu^N_s,\mu_s)\ds + |\nabla\varphi|_{\infty} \int_0^t \|F(\mu^N_s) - F(\mu_s)\|_{L^\infty(B_{R_T}(0))}  \int_{\R^{2d}} \dmu_s(\bx,\bv)\ds. 
\end{split}
\]
The right hand side converges to $0$ as $N\ra \infty$, which implies that
\[
 \int_0^t \int_{\R^{2d}} F(\mu^N_s)(\bx,\bv)\cdot\n_\bv\f(\bx,\bv)  \dmu^N_s(\bx,\bv)\ds \ra \int_0^t \int_{\R^{2d}} F(\mu_s)(\bx,\bv)\cdot\n_\bv\f(\bx,\bv)  \dmu_s(\bx,\bv)\ds
\]
as $N\ra \infty$.
It follows that $\mu$ satisfies Definition \ref{def:ms-soln}. Uniqueness is a direct consequence of the stability \eqref{W1:exp-bound}, which concludes the theorem.  

Note that in the case $\phi$ is fat-tailed and $p\in (2,3)$, $\supp \mu_t^N \subseteq B_{R}(0)$ with $R$ independent of $T$. Thus, in this case there exists a global solution $\mu_t$ with $\supp\mu_t$ is uniformly bounded for all times.
\end{proof}

\section{Propagation of Chaos for Deterministic Model with Fat-tailed Kernels}\label{sec:chaos}
In this section, we prove Theorem \ref{thm:chaos}: the propagation-of-chaos estimate for the deterministic model \eqref{VEfa}:
\begin{equation}\label{VEf}
    \begin{cases}
    \begin{aligned}
     & \p_t f + \bv\cdot \n_{\bx} f + \n_{\bv} \cdot (f F(f)) = 0, \quad f(0) = f_0,\\
     & F(f)(\bx,\bv,t) = \int_{\R^{2d}} \phi(\bx-\by) A(\bw-\bv) f(\by,\bw,t) \dyw.
    \end{aligned}   
    \end{cases}
\end{equation}
Here, we assume without loss of generality that $f_0=f_0(\bx,\bv)$ is a probability distribution.

Recall the Liouville equation \eqref{LEa}:
\begin{equation}\label{LE}
    \begin{cases}
    \begin{aligned}
     & \p_t f^{N} + \sum_{i=1}^N \bv_i \cdot \n_{\bx_i} f^N + \sum_{i=1}^N \n_{\bv_i} \cdot(f^N F_i^N) = 0, \quad f^N(0) = f_0^{\otimes N},\\
     & F_i^N(\bx_1,\ldots,\bx_N,\bv_1,\ldots,\bv_N) = \frac1N \sum_{j=1}^N \phi(\bx_i - \bx_j)A(\bv_j - \bv_i).
    \end{aligned}   
    \end{cases}
\end{equation}
Due to the symmetries of the forces $F_i^N$, the solution $f^N$ will remain symmetric with respect to permutations of pairs $(\bx_i,\bv_i)$ for all time. We denote by 
\[
\Phi^N_t = (X_t, V_t) = \big(\bx_1(t),\ldots,\bx_N(t),\,\bv_1(t),\ldots,\bv_N(t)\big): \R^{2d N} \to \R^{2d N}
\]
the flow-map of \eqref{LE}, in other words, these are solutions to the agent-based system \eqref{NVA}.
Then, $f^N_t$ at any time $t>0$ is a push-forward of the initial distribution by $\Phi^N_t$,
\begin{equation*}\label{}
f^N_t = \Phi_t^N \sharp f_0^{\otimes N}.
\end{equation*}
Now, denote by 
\[
\barPhi_t = (\bar\bx(t),\bar\bv(t)): \R^{2d} \to \R^{2d}
\]
the flow-map of the Vlasov equation \eqref{VEf}, i.e.
\begin{equation*}
\begin{cases}
\begin{aligned}
 \dot{\bar\bx} &= \bar\bv,\\
\dot{\bar\bv} &=  \int_{\R^{2d}} \phi(\bar\bx-\by) A(\bw-\bar\bv) f(\by,\bw,t) \dyw,    
\end{aligned}
\end{cases}
\end{equation*}
and by 
\[
\barPhi^{\otimes N}_t = (\barX_t,\barV_t) = \big(\bar\bx_1(t),\ldots,\bar\bx_N(t), \,\bar\bv_1(t),\ldots, \bar\bv_N(t)\big): \R^{2dN} \to \R^{2dN}
\]
the direct product of $N$ copies of $\barPhi_t$'s.  Thus,
\begin{equation*}\label{}
f_t =  \barPhi_t \sharp f_0, \qquad  f^{\otimes N}_t =  \barPhi_t^{\otimes N} \sharp f_0^{\otimes N}.
\end{equation*}
Denote $\Sigma_N^k$ the set of all different ordered subsets of size $k$ of $\{1,\ldots,N\}$. Then, for any $\s \in \Sigma_N^k$,
\[
 \cW^2_2(f^{(k)}_t,f^{\otimes k}_t)  \leq  \int_{\R^{2dN}} \sum_{i=1}^k | (\bx_{\s(i)},\bv_{\s(i)})-(\bar\bx_{\s(i)},\bar\bv_{\s(i)})|^2 \,f_0^{\otimes N}(X_0,V_0) \dXV.
\]
Summing up over all $\s \in \Sigma_N^k$,
\[
\binom{N}{k} \cW^2_2(f^{(k)}_t,f^{\otimes k}_t)  \leq   \int_{\R^{2dN}} \sum_{\s \in \Sigma_N^k}  \sum_{i=1}^k | (\bx_{\s(i)},\bv_{\s(i)})-(\bar\bx_{\s(i)},\bar\bv_{\s(i)})|^2 \,f_0^{\otimes N}(X_0,V_0) \dXV.
\]
Note that each coordinate is repeated $\binom{N-1}{k-1}$ times on the right hand side of the above inequality. Therefore, we obtain
\begin{align}\label{W2-aux}
 \cW^2_2(f^{(k)}_t,f^{\otimes k}_t) & \leq  \dfrac{k}{N} \int_{\R^{2dN}}  \sum_{i=1}^N | (\bx_{i},\bv_{i})-(\bar\bx_{i},\bar\bv_{i})|^2 \,f_0^{\otimes N}(X_0,V_0) \dXV\notag\\
& =  \dfrac{k}{N}  \int_{\R^{2d N}}  | \Phi_t^N(X_0,V_0) - \barPhi_t^{\otimes N}(X_0,V_0)|^2 \,f_0^{\otimes N}(X_0,V_0) \dXV.
\end{align}
Thus, the proof of \thm{thm:chaos} amounts to establishing the following estimate
\[
  \int_{\R^{2d N}}  | \Phi_t^N(X_0,V_0) - \barPhi_t^{\otimes N}(X_0,V_0)|^2 \,f_0^{\otimes N}(X_0,V_0) \dXV.
\]
Let us break the expression under the integral into potential and kinetic parts,
\begin{align*}
\cP(t) &= \frac12 \int_{\R^{2dN}}  | X_t - \barX_t |^2 \, f_0^{\otimes N} \dXV,\\
 \cK(t) & = \frac12 \int_{\R^{2dN}}  | V_t - \barV_t  |^2\, f_0^{\otimes N}\dXV.
\end{align*}
By the H\"older inequality, we have
\begin{equation}\label{potentialPart}
\begin{split}
 \cP'(t) = \int_{\R^{2dN}}  ( X_t - \barX_t )\cdot(V_t -\barV_t) \, f_0^{\otimes N} \dXV
\leq 2 \sqrt{\cP(t)}\sqrt{\cK(t)}.    
\end{split}
\end{equation}
For the kinetic part, we have
\begin{align*}
 \cK'(t) & =  \int_{\R^{2dN}}\sum_{i=1}^N (\bv_i - \bar\bv_i) \cdot \bigg( \frac1N \sum_{j=1}^N \phi(\bx_i - \bx_j)A(\bv_j - \bv_i)\\
 & \hspace{1in} - \int_{\R^{2d}}\phi(\bar\bx_i-\by) A(\bw-\bar\bv_i) f(\by,\bw,t) \dyw \bigg) \,f_0^{\otimes N} \dXV \\
& =: I_1+I_2+I_3,
\end{align*}
where
\begin{align*}
I_1& = \frac1N\int_{\R^{2dN}} \sum_{i,j=1}^N \big[\phi(\bx_i - \bx_j)- \phi(\bar\bx_i - \bar\bx_j)\big]\,(\bv_i - \bar\bv_i) \cdot A(\bv_j - \bv_i) \, f_0^{\otimes N} \dXV,\\
I_2 & = \frac1N\int_{\R^{2dN}}\sum_{i,j=1}^N \phi(\bar\bx_i - \bar\bx_j)\, (\bv_i - \bar\bv_i) \cdot \big[A(\bv_j-\bv_i) - A(\bar\bv_j-\bar\bv_i)\big] \, f_0^{\otimes N} \dXV,\\
I_3& = \int_{\R^{2dN}}\sum_{i=1}^N (\bv_i - \bar\bv_i) \cdot \bigg( \frac1N \sum_{j=1}^N \phi(\bar\bx_i - \bar\bx_j)A(\bar\bv_j - \bar\bv_i)\\
& \hspace{1in} - \int_{\R^{2d}}\phi(\bar\bx_i-\by) A(\bw-\bar\bv_i) f(\by,\bw,t) \dyw \bigg)\, f_0^{\otimes N}\dXV.
\end{align*}
For $I_1$, we apply the mean value theorem and the smoothness of $\phi$ to obtain 
\begin{align*}
|I_1| & \leq \frac1N\int_{\R^{2dN}} \sum_{i,j=1}^N \|\nabla\phi\|_{\infty}\big[ |\bx_i - \bar\bx_i| + |\bx_j - \bar\bx_j|\big]\,|\bv_i - \bar\bv_i|\,  |\bv_i - \bv_j|^{p-1} \, f_0^{\otimes N} \dXV,\\
& \leq C \max_{i,j=1,\ldots,N} |\bv_i - \bv_j|^{p-1} \sqrt{\cK(t)}\, \Big(2 \int_{\R^{2d N}} \sum_{i=1}^N |\bx_i -\bar\bx_i|^2 \, f_0^{\otimes N} \dXV\Big)^{1/2}\\
& \leq C\,\cV_\Omega^{p-1}(t)\, \sqrt{\cK(t)} \sqrt{\cP(t)}.
\end{align*}
Next, for $I_2$, since $\phi$ is even, and $A$ is odd, symmetrizing $i$ and $j$ yields
\begin{align*}
I_2 & = \frac{1}{2N}\int_{\R^{2dN}}\sum_{i,j=1}^N \phi(\bar\bx_i - \bar\bx_j)\,[(\bv_i - \bar\bv_i) - (\bv_j - \bar\bv_j)] \cdot \big[A(\bv_j-\bv_i) - A(\bar\bv_j-\bar\bv_i)\big] \, f_0^{\otimes N} \dXV\leq0,
\end{align*}
where we have used \eqref{eq:eleineq2} with $\bx=\bv_j-\bv_i$ and $\by = \bar\bv_j-\bar\bv_i$.
Finally, for $I_3$, applying the \HI\ we have
\[
I_3 \leq \Big(\int_{\R^{2dN}} \sum_{i=1}^N |\bv_i- \bar \bv_i|^2f_0^{\otimes N} \dXV \Big)^{\frac12} J^{\frac12} = \sqrt{\cK(t)}\,\sqrt{J},
\]
where we further estimate
\begin{align*}
J & =\int_{\R^{2dN}} \sum_{i=1}^N \Big|\frac1N \sum_{j=1}^N \phi(\bar\bx_i - \bar\bx_j)A(\bar\bv_j - \bar\bv_i) - \int_{\R^{2d}}\phi(\bar\bx_i-\by) A(\bw-\bar\bv_i) f(\by,\bw,t) \dyw  \Big|^2 f_0^{\otimes N} \dXV\\
& = \sum_{i=1}^N \int_{\R^{2dN}} \Big|\frac1N \sum_{j=1}^N \phi(\bar\bx_i - \bar\bx_j)A(\bar\bv_j - \bar\bv_i) - \int_{\R^{2d}}\phi(\bar\bx_i-\by) A(\bw-\bar\bv_i) f(\by,\bw,t) \dyw  \Big|^2 f_t^{\otimes N} \dd\bar{X}_t \dd\bar{V}_t\\
& \leq N \, \Big( \frac4N \sup_{(\bar\bx',\bar\bv'),(\bar\bx'',\bar\bv'')\in \supp f_t} |\phi(\bar\bx'-\bar\bx'')A(\bar\bv'-\bar\bv'')|^2\Big)
\leq 4\|\phi\|_\infty^2 \cV_\Omega^{2(p-1)}(t).
\end{align*}
For the penultimate inequality, we have used the estimate in \cite[Lemma 3.3]{natalini2021mean} to control each integrand. Therefore, we obtain
\begin{equation*}
I_3 \le C \cV^{p-1}_{\Omega}(t)\sqrt{\cK(t)}.
\end{equation*}
Summing up the above estimates for $I_1, I_2$ and $I_3$, we deduce that
\begin{equation}\label{kinPart}
 \cK'(t) \leq C \cV^{p-1}_{\Omega}(t) \sqrt{\cK(t)} \big(1 + \sqrt{\cP(t)}\big).
\end{equation}
Set
\[
 a(t) := 1+ \sqrt{\cP(t)},\qquad b(t) := \sqrt{\cK(t)}.
\]
Then, we deduce from  \eqref{potentialPart} and \eqref{kinPart} that the dynamics of $(a,b)$ satisfies the paired inequalities
\begin{equation}\label{auxSys}
\begin{cases}
    \begin{aligned}
     a'(t) &\leq b(t),\quad\ a(0) = 1,\\ 
     b'(t) &\leq C \cV^{p-1}_{\Omega}(t)\, a(t),\quad b(0) = 0.
    \end{aligned}
\end{cases}
\end{equation}
The term $\cV^{p-1}_{\Omega}(t)$ can be further estimated by Theorem \ref{thm:kin-flocking}. Hence, the system \eqref{auxSys} has the form \eqref{ab-aux-sys} with $g(t)\equiv0$ (or \eqref{ab-aux-sys2} for the borderline case when $p=3$). We apply Lemma \ref{lem:ab-aux} (and Lemma \ref{lem:ab-aux2} for $p=3$) and obtain the bounds analogous to those in Lemma \ref{claim:grad-est}:
 \begin{enumerate}[label=\textnormal{(\roman*).}]
    \item if $p\in (2,3)$ then
    \[
	\cP(t) \lesssim \lan t\ran^2,\quad
    \cK(t) \lesssim 1; 
    \]
   \item if $p>3$ then 
    \[
    \cP(t) \lesssim e^{2\bar C\lan t\ran^{\gamma}}, \quad  
	\cK(t) \lesssim \lan t\ran^{-2(1-\gamma)} e^{2\bar C\lan t\ran^{\gamma}},\quad \gamma = \tfrac{(1+\alpha)(p-3)}{2(p-\alpha-2)};
    \]
   \item if $p=3$ then
    \[
	\cP(t) \lesssim \, e^{2\bar C(\log \langle t\rangle)^{\frac{1}{1-\alpha}}}, \quad 
	\cK(t) \lesssim \, \langle t\rangle^{-2}(\log \langle t\rangle)^{\frac{2\alpha}{1-\alpha}}\,e^{2\bar C (\log \langle t\rangle)^{\frac{1}{1-\alpha}}}.
    \] 
\end{enumerate}  
Note that in the case of $p\in (2,3)$, thanks to the flocking estimate in Theorem \ref{thm:kin-flocking} (i), we also have $\cP(t)\leq CN$. Thus, in this case
\[
\cP(t) \lesssim \min\{ N, t^{2}\}.
\]
Plugging all the bounds above into \eqref{W2-aux}, we conclude with the propagation-of-chaos estimates:
\begin{enumerate}[label=\textnormal{(\roman*).}]
    \item if $p\in (2,3)$ then
    \[
    \cW_2(f^{(k)}_t,f^{\otimes k}_t) \lesssim \sqrt{k}\min\left\{1,\frac{t}{\sqrt{N}}\right\};
    \]
    \item if $p > 3$ then
     \[
      \cW_2(f^{(k)}_t,f^{\otimes k}_t) \leq C\sqrt{\frac{k}{N}}e^{\bar C \langle t\rangle^{\frac{(1+\alpha)(p-3)}{2(p-\alpha-2)}}};
     \]
      \item if $p=3$ then
     \[
\cW_2(f^{(k)}_t,f^{\otimes k}_t) \leq C\sqrt{\frac{k}{N}}(\log \langle t\rangle)^{\frac{\alpha}{1-\alpha}}e^{\bar C (\log\langle t\rangle)^{\frac{1}{1-\alpha}}}.
     \]
\end{enumerate}
This finishes the proof of Theorem \ref{thm:chaos}.

\section{Mean-Field Limit for Stochastic Model}\label{sec:StochasticMFL}
In this section we establish similar results in the stochastic setting and prove Theorem~\ref{thm:StochasticMFL}.

Let $(\bx_i^0,\bv_i^0)$ be $N$ independent identically distributed (i.i.d.) random variables corresponding to law  $f_0$, and let $(\bx_i,\bv_i)$ be the solution to \eqref{SNVA}.  Both conclusions of the theorem come from a single source -- comparing solutions of \eqref{SNVA} with $N$ identically distributed and independent versions of the stochastic characteristics of the Fokker-Planck equation \eqref{FPE}. So, let $(\bar\bx_i, \bar\bv_i), i = 1, \ldots, N$ solve:
\begin{equation}\label{FPEchar}
\begin{cases}
\begin{aligned}
 \dd{\bar\bx_i} &= \bar\bv_i \dt,\\
\dd{\bar\bv_i} &=  \Big(\int_{\dom} \phi(\by- \bar\bx_i) A(\bw-\bar\bv_i) f(\by,\bw,t) \dd\by \dd\bw\Big) \dd t + \sqrt{2h(s_\rho(\bar\bv_i))}\dd W_i, \end{aligned}
\end{cases}
\end{equation}
where
\[
s_\rho(\bar\bx_i) := \phi \ast\rho ( \bar\bx_i), \quad \rho(\bx) = \int_{\R^d}f(\bx,\bv)\dd \bv,
\]
and $W_i$ are independent Brownian motions. We also let initial conditions $(\bx_i^0,\bv_i^0)$ have a common law $f_0$ and be independent. As a result,  $(\bar\bx_i,\bar\bv_i)$ are i.i.d.'s with the common law  $f$.

We have
\begin{align*}
& \E \Big|  \frac1N \sum_{i=1}^N \varphi(\bx_i(t),\bv_i(t)) - \int_{\dom} \varphi(\bx,\bv) f(\bx,\bv,t) \dd\bx\dd\bv \Big|^2 \\ 
 \leq & \,\E  \Big|  \frac1N \sum_{i=1}^N \varphi(\bx_i(t),\bv_i(t)) - \frac1N \sum_{i=1}^N \varphi(\bar\bx_i(t),\bar\bv_i(t))  \Big|^2  \\
& +  \underbrace{\E  \Big|  \frac1N \sum_{i=1}^N \varphi(\bar\bx_i(t),\bar\bv_i(t)) - \int_{\dom} \varphi(\bx,\bv) f(\bx,\bv,t) \dd\bx\dd\bv \Big|^2}_{=:J} .
\end{align*}
Denoting
\[
\mu_{\varphi}(t):=\int_{\dom}\varphi(\bx,\bv)\,f(\bx,\bv,t)\dd\bx \dd\bv
\]
and $\varphi_i:=\varphi(\bar\bx_i(t),\bar\bv_i(t))$. Then $\varphi_i$'s are i.i.d.\ with mean $\mu_{\varphi}(t)$, and
\[
J=\E\Big|\frac1N\sum_{i=1}^N\varphi_i-\mu_\varphi(t)\Big|^2.
\]
Expanding the square and using independence,
\[
\begin{split}
J
&=\E\Bigg[\frac{1}{N^2}\sum_{i=1}^N(\varphi_i-\mu_{\varphi}(t))^2
+\frac{1}{N^2}\sum_{i\neq j}(\varphi_i-\mu_{\varphi}(t))(\varphi_j-\mu_{\varphi}(t))\Bigg]\\
&=\frac{1}{N^2}\sum_{i=1}^N \E\big[(\varphi_i-\mu_{\varphi}(t))^2\big]
+\frac{1}{N^2}\sum_{i\neq j}\E[\varphi_i-\mu_{\varphi}(t)]\E[\varphi_j-\mu_{\varphi}(t)] = \frac{1}{N}\Var_f(\varphi),
\end{split}
\]
where
\[
\Var_f(\varphi):=\int_{\dom}\varphi(\bx,\bv)^2f(\bx,\bv,t)\dd\bx \dd\bv
-\Big(\int_{\dom}\varphi(\bx,\bv)f(\bx,\bv,t)\dd\bx \dd\bv\Big)^2.
\]
By the symmetry, one has
\[
\E \Big|  \frac1N \sum_{i=1}^N \varphi(\bx_i(t),\bv_i(t)) - \frac1N \sum_{i=1}^N \varphi(\bar\bx_i(t),\bar\bv_i(t))  \Big|^2  \leq  \Lip(\varphi)^2 \E[ |\bx_i - \bar\bx_i|^2 + |\bv_i - \bar\bv_i|^2].
\]

At the same time, following \cite{bahouri2011fourier}, 
\[
W_2^2(f^{(k)}, f^{\otimes k}) \leq \E\left[ \sum_{i=1}^k |x_i - \barx_i|^2 + |v_i - \barv_i|^2 \right] = k  \E[ |\bx_i - \bar\bx_i|^2 + |\bv_i - \bar\bv_i|^2].
\]
Note that in the the above the right hand sides are independent of particular $i$. 

Thus, denoting
\[
\cE(t) := \E[ |\bx_i - \bar\bx_i|^2 + |\bv_i - \bar\bv_i|^2] =: \cE_{\bx}(t) + \cE_{\bv}(t).
\]
the proof of the theorem reduces to establishing the bound on expected divergence of characteristics
\[
\cE(t) \leq  \frac{C}{N^{e^{-Ct}}}.
\]

\begin{proof}[Proof of \thm{thm:StochasticMFL}]
Compute the time derivative of  $\cE_\bx$, we have
\begin{equation}\label{Ex}
  \ddt \cE_\bx(t) = 2  \E[(\bx_i - \bar\bx_i) \cdot (\bv_i - \bar\bv_i) ] \leq \cE(t).  
\end{equation}
By It\^o's formula and systems \eqref{SNVA}, \eqref{FPEchar}, the time derivative of $\cE_\bv(t)$ is computed as
\begin{align*}
\dfrac{1}{2}\cE_\bv'(t) = &\, \E \Big[ (\bv_i - \bar\bv_i)\cdot\Big(\dfrac{1}{N}\sum\limits_{j= 1}^N \phi(\bx_j - \bx_i)A(\bv_j- \bv_i) - \int_{\dom} \hspace{-.1in}\phi(\by- \bar\bx_i) A(\bw-\bar\bv_i) f(\by,\bw,t) \dd\by \dd\bw\Big)\Big]\\
&\quad +\E\left[ (\bv_i - \bar\bv_i)\cdot \sqrt{2}\big(\sqrt{h(s_i)}-\sqrt{h(s_\rho(\bar\bx_i)}\big)\,\dd W_i\right]
+ d\,\E\Big[\Big(\sqrt{h(s_i)}-\sqrt{h(s_\rho(\bar\bx_i))}\Big)^2\Big] \\
= &\, \dfrac{1}{N}\sum\limits_{j= 1}^N \E \Big[\phi(\bx_j - \bx_i)(\bv_i - \bar\bv_i)\cdot\Big(A(\bv_j- \bv_i) - A(\bar\bv_j-\bar\bv_i) \Big) \Big]\\
 & +\dfrac{1}{N}\sum\limits_{j= 1}^N \E \Big[(\bv_i - \bar\bv_i)\cdot A(\bar\bv_j-\bar\bv_i) \Big(\phi(\bx_j - \bx_i) - \phi(\bar\bx_j - \bar\bx_i) \Big)\Big]\\
& +\dfrac{1}{N}\sum\limits_{j= 1}^N \E \Big[(\bv_i - \bar\bv_i)\cdot \Big(\phi(\bar\bx_j - \bar\bx_i)A(\bar\bv_j-\bar\bv_i) -  \int_{\dom} \hspace{-.2in}\phi(\by- \bar\bx_i) A(\bw-\bar\bv_i) f(\by,\bw,t) \dd\by \dd\bw\Big)\Big]\\
&+\E\left[ (\bv_i - \bar\bv_i)\cdot \sqrt{2}\big(\sqrt{h(s_i)}-\sqrt{h(s_\rho(\bar\bx_i)}\big)\,\dd W_i\right]\\
&+ d\,\E\Big[\Big(\sqrt{h(s_i)}-\sqrt{h(s_\rho(\bar\bx_i))}\Big)^2\Big] \\
:= &\, J_1+J_2+J_3+J_4+J_5.
\end{align*}
Because agents are independent and identically distributed, $\phi$ is even (radial), and $A$ is odd, $J_1$ can be rewritten as
\begin{align}
J_1 & = \dfrac{1}{2N^2}\sum\limits_{i,j= 1}^N \E \Big[\phi(\bx_j - \bx_i)\big((\bv_i - \bar\bv_i) - (\bv_j - \bar\bv_j)\big)\cdot\Big(A(\bv_j- \bv_i) - A(\bar\bv_j-\bar\bv_i) \Big) \Big]\notag\\
& = -\dfrac{1}{2N^2}\sum\limits_{i,j= 1}^N \E \Big[\phi(\bx_j - \bx_i)\big((\bv_j - \bv_i) - (\bar\bv_j - \bar\bv_i)\big)\cdot\Big(A(\bv_j- \bv_i) - A(\bar\bv_j-\bar\bv_i) \Big) \Big] \leq 0. \label{J1}  
\end{align}
For $J_2$, we have
\[
\begin{split}
    & J_2 = \dfrac{1}{N}\sum\limits_{j= 1}^N \E \Big[(\bv_i - \bar\bv_i)\cdot \big|\bar\bv_j-\bar\bv_i\big|^{p-2}(\bar\bv_j-\bar\bv_i) \Big(\phi(\bx_j - \bx_i) - \phi(\bar\bx_j - \bar\bx_i) \Big)\Big]\\
    &\quad\leq C \E\Big[ \big|\bv_i - \bar\bv_i\big| \big|\bar\bv_j-\bar\bv_i \big|^{p-1}\min\big\{1,|(\bx_j - \bx_i) -(\bar\bx_j - \bar\bx_i) |\big\}\Big] 
\end{split}
\]
since $\phi$ is bounded above and smooth. Let us fix a pair $(i,j)$ and define the random variable
\[
\bz := \big|\bv_i - \bar\bv_i \big| \,\big|\bar\bv_j-\bar\bv_i \big|^{p-1}\min\big\{1,|(\bx_j - \bx_i) -(\bar\bx_j - \bar\bx_i) |\big\}.
\]
For any given positive constant $R$, letting
\[
\cR := \big\{ |\bar\bv_i|\leq R, |\bar\bv_j| \leq R\big\}.
\]
Then, 
\[
    \E[\bz]  = \E[\one_{\cR}\bz] + \E[\one_{\cR^c}\bz]. 
\]
For the first term, we have
\[
\begin{split}
    \E[\one_{\cR}\bz]&\lesssim R^{p-1}\E\Big[\big|\bv_i - \bar\bv_i \big|\,\big |(\bx_j - \bx_i) -(\bar\bx_j - \bar\bx_i) \big|\Big] \\
    &\lesssim R^{p-1} \cE(t) \quad \text{ (by the Cauchy-Schwarz inequality)}.
\end{split}
\]
Using Young's inequality, the second term is estimated by
\[
\begin{split} 
    \E[\one_{\cR^c}\bz]& \lesssim \E\Big[\one_{\cR^c}\big|\bv_i - \bar\bv_i \big|^2\Big] + \E \Big[\one_{\cR^c} \big|\bar\bv_j-\bar\bv_i \big|^{2(p-1)}\Big]\\
    &\lesssim  \E\Big[\big|\bv_i - \bar\bv_i\big|^2\Big] + \big(\E \big[\one_{\cR^c}\big]\big)^{1/2}\Big(\E\big[ \big|\bar\bv_j-\bar\bv_i \big|^{4(p-1)}\big]\Big)^{1/2}\\
    & \lesssim \cE(t) + \Big(\E \big[\one_{\{|\bar\bv_i| > R\}}\big] + \E \big[\one_{\{|\bar\bv_j| > R\}}\big]\Big)^{1/2}\Big(\E\big[ \big|\bar\bv_i \big|^{4(p-1)}\big]\Big)^{1/2}.
\end{split}  
\]
By Markov's inequality,
\[
\begin{split}
    \E \big[\one_{\{|\bar\bv_i| > R\}}\big] & \leq \P\big(e^{c_0|\bar\bv_i|} > e^{c_0R^{p-1}}\big) \leq \frac{1}{e^{c_0R^{p-1}}}  \E \big[e^{c_0|\bar\bv_i|}\big].
\end{split}
\]
Since $\E \big[e^{c_0|\bar\bv_i|}\big]< +\infty$, one has
\[
\E \big[\one_{\{|\bar\bv_i| > R\}}\big]  \lesssim e^{-c_0R^{p-1}}, \quad \E\big[ \big|\bar\bv_i \big|^{4(p-1)}\big] \leq \E \big[e^{c_0|\bar\bv_i|}\big]< +\infty. 
\]
Hence,
\[
 \E[\one_{\cR^c}\bz] \lesssim \cE(t) + e^{-\frac{c_0R^{p-1}}{2}}.
\]
Therefore,
\begin{equation}\label{J2}
    J_2 \lesssim (1+R^{p-1})\cE(t) + e^{-\frac{c_0 R^{p-1}}{2}}.
\end{equation}
For $J_3$, we have
\[
\begin{split}
J_3 &= \dfrac{1}{N}\sum\limits_{j\ne i}^N \E \Big[(\bv_i - \bar\bv_i)\cdot \Big(\phi(\bar\bx_j - \bar\bx_i)A(\bar\bv_j-\bar\bv_i) -  \int_{\dom} \phi(\by- \bar\bx_i) A(\bw-\bar\bv_i) f(\by,\bw,t) \dd\by \dd\bw\Big)\Big] \\
&\quad - \dfrac{1}{N} \E \Big[(\bv_i - \bar\bv_i)\cdot\int_{\dom} \phi(\by- \bar\bx_i) A(\bw-\bar\bv_i) f(\by,\bw,t) \dd\by \dd\bw\Big]:= J_{31}+J_{32}.
\end{split}
\]
By Cauchy-Schwarz inequality,
\begin{equation}\label{J32}
 J_{32}  \leq \frac{1}{N} \Big(\E\Big[\big|\bv_i - \bar\bv_i\big|^2\Big]\Big)^{1/2} \, \Big(\E\Big[\int_{\dom} \phi(\by- \bar\bx_i) A(\bw-\bar\bv_i) f(\by,\bw,t) \dd\by \dd\bw\Big]\Big)^{1/2}
\lesssim \frac{\sqrt{\cE(t)}}{N}.   
\end{equation}
Here we use the fact that
\[
\begin{split}
  &\E\Big[\int_{\dom} \phi(\by- \bar\bx_i) A(\bw-\bar\bv_i) f(\by,\bw,t) \dd\by \dd\bw\Big]\\
   & = \int_{(\dom)^2}   \phi(\by- \bx) A(\bw-\bv) f(\by,\bw,t)  f(\bx,\bv,t)\dd\by \dd\bw \dd\bx\dd\bv\\
  & \leq \|\phi\|_{\infty} \int_{(\dom)^2}|\bw -\bv|^{p-1}f(\by,\bw,t)  f(\bx,\bv,t)\dd\by \dd\bw \dd\bx\dd\bv\\
  & \lesssim  \int_{(\dom)^2}\big(|\bw|^{p-1} +|\bv|^{p-1}\big)f(\by,\bw,t)  f(\bx,\bv,t)\dd\by \dd\bw \dd\bx\dd\bv\\
  & \lesssim  \int_{\dom}|\bv|^{p-1}  f(\bx,\bv,t) \dd\bx\dd\bv <+\infty \quad \text{ (due to \eqref{moment})}.
\end{split}
\]
Let us consider $J_{31}$. For each $ j \neq i $, define
\[
\xi_j := \phi(\bar\bx_j - \bar\bx_i)A(\bar\bv_j-\bar\bv_i) - \int_{\dom} \phi(\by- \bar\bx_i) A(\bw-\bar\bv_i) f(\by,\bw,t) \dd\by \dd\bw.
\]
Then,
\[
J_{31} = \dfrac{1}{N} \E \Big[ (\bv_i - \bar{\bv}_i) \cdot \Big( \sum_{j \neq i}^N \xi_j \Big) \Big].
\]
Applying the Cauchy-Schwarz inequality,
\begin{align}
J_{31} & \leq \dfrac{1}{N} \E \Big[ |\bv_i - \bar{\bv}_i| \, \Big | \sum_{j \neq i}^N \xi_j \Big | \Big ]\notag \\
& \leq \dfrac{1}{N} \left( \E\big[|\bv_i - \bar{\bv}_i|^2\big] \right)^{1/2} \Big( \E\Big[ \Big | \sum_{j \neq i}^N \xi_j \Big |^2 \Big ] \Big)^{1/2}\notag\\
& \leq \dfrac{\sqrt{\cE_v(t)} }{N} \Big( \E\Big[ \Big | \sum_{j \neq i}^N \xi_j \Big |^2 \Big ] \Big)^{1/2}.   \label{J31:aux1}
\end{align}
We now examine the variance
\[
\E\Big[ \Big | \sum_{j \neq i}^N \xi_j \Big |^2 \Big ] = \sum_{j \neq i}^N \E[|\xi_j|^2] + \sum_{\substack{j, k \neq i \\ j \neq k}} \E[\xi_j \cdot \xi_k].
\]
Since $\{(\bar{\bx}_j, \bar{\bv}_j)\}_{j \neq i}$ are i.i.d. with law $f$, and $\xi_j$ depends on $(\bar{\bx}_j, \bar{\bv}_j)$ and $(\bar{\bx}_i, \bar{\bv}_i)$, we have for $j \neq k$,
\[
\E[\xi_j \cdot \xi_k] = \E\Big[ \E\big[\xi_j \cdot \xi_k \mid (\bar{\bx}_i, \bar{\bv}_i)\big] \Big ].
\]
Conditional on $(\bar{\bx}_i, \bar{\bv}_i)$, $\xi_j$ and $\xi_k$ are independent for $j \neq k$ and have zero mean
\[
\E[\xi_j \mid (\bar{\bx}_i, \bar{\bv}_i)] = 0.
\]
Therefore,
\[
\E[\xi_j \cdot \xi_k \mid (\bar{\bx}_i, \bar{\bv}_i)] = \E[\xi_j \mid (\bar{\bx}_i, \bar{\bv}_i)] \cdot \E[\xi_k \mid (\bar{\bx}_i, \bar{\bv}_i)] = 0.
\]
Hence, $\E[\xi_j \cdot \xi_k] = 0$ for $j \neq k$. Thus, the cross-terms vanish, and we have
\[
\E\Big[ \Big | \sum_{j \neq i}^N \xi_j \Big |^2 \Big ] = \sum_{j \neq i}^N \E[|\xi_j|^2].
\]
By exchangeability, $\E[|\xi_j|^2]$ is the same for all $j \neq i$. Fix $j_0 \ne i$, then
\begin{equation}\label{J31:aux2}
 \E\Big[ \Big | \sum_{j \neq i}^N \xi_j \Big |^2 \Big ] = (N-1) \E[|\xi_{j_0}|^2].   
\end{equation}
Since $ \E[\xi_{j_0} \mid (\bar{\bx}_i, \bar{\bv}_i)] = 0 $, one has
\[
 \E[|\xi_{j_0}|^2] = \E\Big[ \Var\Big( \phi(\bar\bx_{j_0} - \bar\bx_i)A(\bar\bv_{j_0}-\bar\bv_i) \mid (\bar{\bx}_i, \bar{\bv}_i) \Big) \Big ].
\]
Using the boundedness of $ \phi $ and the growth of $ A $, one gets
\[
|\phi(\bar\bx_{j_0} - \bar\bx_i)A(\bar\bv_{j_0}-\bar\bv_i)| \leq \|\phi\|_{\infty} |\bar{\bv}_{j_0} - \bar{\bv}_i|^{p-1} \lesssim (|\bar{\bv}_{j_0}|^{p-1} + |\bar{\bv}_i|^{p-1}).
\]
Therefore, the conditional variance is bounded by the conditional second moment
\[
\Var\Big( \phi(\bar\bx_{j_0} - \bar\bx_i)A(\bar\bv_{j_0}-\bar\bv_i) \mid (\bar{\bx}_i, \bar{\bv}_i) \Big) \leq \E\Big[ |\phi(\bar\bx_{j_0}- \bar\bx_i)A(\bar\bv_{j_0}-\bar\bv_i)|^2 \mid (\bar{\bx}_i, \bar{\bv}_i) \Big ].
\]
Thus,
\[
\ \E[|\xi_{j_0}|^2]  \leq \E\Big[ \big|\phi(\bar\bx_{j_0} - \bar\bx_i)A(\bar\bv_{j_0}-\bar\bv_i)\big|^2 \Big ] \lesssim \E\Big[ (|\bar{\bv}_{j_0}|^{p-1} + |\bar{\bv}_i|^{p-1})^2 \Big ] < +\infty.
\]
Combining this inequality with \eqref{J31:aux1} and \eqref{J31:aux2}, we have
\begin{equation*}
    J_{31} \leq \frac{C\sqrt{N-1}\sqrt{\cE_{\bv}(t)}}{N} \leq  \frac{C\sqrt{\cE_{\bv}(t)}}{\sqrt{N}} \leq \cE_{\bv}(t) + \frac{C}{N}.
\end{equation*}
This inequality and \eqref{J32} imply that 
\begin{equation}\label{J3}
    J_{3} \leq \cE(t) + \frac{C}{N}.
\end{equation}
Because 
\[
M_t:= \int_0^t\big(\bv_i(\tau) - \bar\bv_i(\tau)\big)\cdot \sqrt{2}\left(\sqrt{h(s_i(\tau))}-\sqrt{h(s_\rho(\bar\bx_i(\tau))}\right)\,\dd W_i(\tau)
\]
is a martingale, so $\E[M_t] = \E[M_0] = 0$. Thus,
\begin{equation}\label{J4}
   J_4 = 0.
\end{equation}
For $J_5$, since $h$ is Lipschitz and due to \eqref{pos} we have
\[
    J_5 = d\,\E\Big[\Big|\frac{h(s_i)-h(s_\rho(\bar\bx_i))}{\sqrt{h(s_i)}+\sqrt{h(s_\rho(\bar\bx_i))}}\Big|^2\Big]
    \leq  d\,\Lip(h)^2\E\Big[\big|s_i - s_\rho(\bar\bx_i)\big|^2\Big]
\]
We have
\begin{align*}
|s_i - s_\rho(\bar\bx_i)|
&= \Big| \frac{1}{N} \sum_{j=1}^N \phi(\bx_i - \bx_j) - \int_{\dom} \phi(\bar\bx_i - \by) f(\by, \bw,t)\, \dd\by\, \dd\bw \Big| \\
&\le \Big|  \frac{1}{N} \sum_{j=1}^N \big(\phi(\bx_i - \bx_j) - \phi(\bar\bx_i - \bar\bx_j)\big)\Big| \\
&\quad + \Big| \frac{1}{N} \sum_{j=1}^N \phi(\bar\bx_i - \bar\bx_j) - \int_{\dom} \phi(\bar\bx_i - \by) f(\by, \bw,t)\, \dd\by\, \dd\bw \Big|.
\end{align*}
Thus, 
\[
\begin{split}
    \E\Big[\big|s_i - s_\rho(\bar\bx_i)\big|^2\Big] \leq &\, 2\|\nabla\phi\|_\infty^2\E \Big[\Big| \frac{1}{N}\sum_{j=1}^N\big|(\bx_i -\bar\bx_i) - (\bx_j -\bar\bx_j)\big|\Big|^2 \Big]\\
   & +2\E \Big[\Big| \frac{1}{N} \sum_{j=1}^N \phi(\bar\bx_i - \bar\bx_j) - \int_{\dom} \phi(\bar\bx_i - \by) f(\by, \bw,t)\, \dd\by\, \dd\bw \Big|^2\Big]\\
   =: & \,J_{51}+J_{52}.
\end{split}
\]
We get
\begin{equation}\label{J51}
 J_{51}\leq C\E\Big[|\bx_i - \bar\bx_i|^2\Big] \leq C\cE_x(t) \quad \text{ (due to i.i.d. property)}.   
\end{equation}
For $J_{52}$, let us define for each $j\ne i$,
\[
\zeta_j: = \phi(\bar\bx_i - \bar\bx_j) - \int_{\dom} \phi(\bar\bx_i - \by) f(\by, \bw,t)\, \dd\by\, \dd\bw,
\]
then
\begin{equation}\label{J52:aux1}
 J_{52}\leq \frac{2\|\phi\|_{\infty}}{N^2} + \frac{2}{N^2} \sum_{j\ne i}\E\left[|\zeta_j|^2\right] + \frac{2}{N^2} \sum_{\substack{j, k \neq i \\ j \neq k}}  \E\left[\zeta_j\cdot\zeta_k\right]   
\end{equation}
Since $\{\bar{\bx}_j, \bar{\bv}_j)\}_{j \neq i}$ are i.i.d. with law $f$, and $\zeta_j$ depends on $\bar{\bx}_j$ and $\bar{\bx}_i$, for $j \neq k$,
\[
\E[\zeta_j \cdot \zeta_k] = \E\Big[ \E\big[\zeta_j \cdot \zeta_k \mid \bar{\bx}_i \big] \Big ].
\]
Conditional on $\bar{\bx}_i$, $\zeta_j$ and $\zeta_k$ are independent for $j \neq k$ and have zero mean, hence,
\[
 \E\Big[ \E\big[\zeta_j \cdot \zeta_k \mid \bar{\bx}_i \big] \Big ] =  \E\Big[ \E\big[\zeta_j \mid \bar{\bx}_i \big]   \E\big[ \zeta_k \mid \bar{\bx}_i \big] \Big ] = 0, 
\]
By exchangeability, fix $j_0 \ne i$, we have
\begin{equation}\label{J52:aux2}
 \sum_{j\ne i}\E\left[|\zeta_j|^2\right]  = (N-1)\E\left[|\zeta_{j_0}|^2\right].   
\end{equation}
Since $ \E[\zeta_{j_0} \mid \bar{\bx}_i] = 0 $, 
\[
 \E[|\xi_{j_0}|^2] = \E\Big[\E[|\xi_{j_0}|^2 \mid \bar\bx_i] \Big] = \E\Big[ \Var\big( \phi(\bar\bx_{j_0} - \bar\bx_i)\mid \bar{\bx}_i \big) \Big ].
\]
As $\phi$ is bounded, 
\[
|\phi(\bar\bx_{j_0} - \bar\bx_i)| \leq \|\phi\|_{\infty}.
\]
Therefore, the conditional variance is bounded by the conditional second moment
\[
\Var\Big( \phi(\bar\bx_{j_0} - \bar\bx_i) \mid \bar{\bx}_i \Big) \leq \E\Big[ |\phi(\bar\bx_{j_0}- \bar\bx_i)|^2 \mid \bar{\bx}_i \Big ].
\]
Thus,
\[
\ \E[|\zeta_{j_0}|^2]  \leq \E\Big[ \big|\phi(\bar\bx_{j_0} - \bar\bx_i)\big|^2 \Big ] < +\infty.
\]
This inequality together with \eqref{J52:aux1} and \eqref{J52:aux2} implies
\begin{equation*}
    J_{52} \leq \frac{C}{N}.
\end{equation*}
Combining this with \eqref{J51}, we yield
\begin{equation}\label{J5}
    J_5 \leq \frac{C}{N} + C\cE_x(t).
\end{equation}
By estimates \eqref{Ex}, \eqref{J1}, \eqref{J2}, \eqref{J3}, \eqref{J4}, and \eqref{J5}  we obtain
\begin{equation*}
    \ddt \cE(t) \leq C(1+R^{p-1})\cE(t) + e^{-\frac{c_0\,R^{p-1}}{2}} + \frac{C}{N} \quad \forall t\in [0,T].
\end{equation*}
It implies that for any $r>0$, there exists constant $C_0$ such that
\begin{equation}\label{dE}
\ddt \cE(t) \leq C_0(1+r)\cE(t) + e^{-r} + \frac{C_0}{N} \quad \forall t\in [0,T].    
\end{equation}
By Gr\"onwall's lemma,
\[
\cE(t) \leq \cE(0) e^{C_0(1+r)t} + \int_0^t e^{C_0(1+r)(t-s)}\Big(e^{-r} + \frac{C_0}{N}\Big)\dd s \quad \forall t\in [0,T].
\]
Thus, for all $N\geq 1$, 
\[
    \cE(t) \leq (e^{-r} +C_0)(e^{C_0(1+r)T}-1):= C_1 \quad \forall t\in [0,T].
\]
We now define 
\[
a(t) := \frac{\cE(t)}{e C_1},
\]
which ensures that $a(t) \leq e^{-1}$. This bound implies the useful inequality 
\[
1 - \ln a(t) \leq -2 \ln a(t).
\]
For any $t$ with $a(t) > 0$, we choose $r = -\ln a(t) > 0$. Substituting this into \eqref{dE} yields
\begin{equation}\label{e:a}
a'(t) \leq C_0 \bigl(1 - \ln a(t)\bigr)\, a(t) + \frac{e C_0 C_1}{N}
\leq -C_2\, a(t)\, \ln a(t) + \frac{C_2}{N},
\end{equation}
with $C_2 = \max \{2C_0, e C_0 C_1\}$.
The same inequality holds trivially when $a(t) = 0$ by the standard convention $0 \ln 0 = 0 $. Thus, \eqref{e:a} is valid for all $t \in [0, T]$.
To simplify the analysis, we rescale time by defining
\[
u(t) := a(C_2 t).
\]
Then $ u(0) = 0$, and \eqref{e:a} becomes
\[
u'(t) \leq -u(t) \ln u(t) + \frac{1}{N}, \qquad \text{for } t \in [0, T/C_2].
\]
We now introduce a time-dependent scaling factor to absorb the $1/N$ term. Let $b(t)$ be a function to be determined and define
\[
v(t) := u(t)\, N^{b(t)}.
\]
A computation shows that $v(t)$ satisfies $v(0) = 0$ and
\[
v'(t) \leq -v(t) \ln v(t) + N^{b(t)-1} + v(t)\, \ln N \cdot \bigl(b'(t) + b(t)\bigr).
\]
The optimal choice is $b(t) = e^{-t}$, which is bounded by 1. With this selection, the inequality reduces to
\[
v'(t) \leq -v(t) \ln v(t) + 1 \leq 1 + e^{-1},
\]
where the last bound follows from maximizing $-z \ln z $ for $z>0$.
Integrating this differential inequality from $0$ to $t$ gives
\[
v(t) \leq \left(1 + \frac{1}{e}\right) t   \leq \left(1 + \frac{1}{e}\right)\frac{T}{C_2}  \quad \text{for } t\in [0, T/C_2].
\]
Reverting to the original variables, we conclude that for all $t \in [0, T]$, 
\[
\mathcal{E}(t) = \mathbb{E}\left[\,|\mathbf{x}_i - \bar{\mathbf{x}}_i|^2 + |\mathbf{v}_i - \bar{\mathbf{v}}_i|^2 \,\right] \leq C_3\, N^{-e^{-C_4 t}},
\]
with $C_3 = \left(1 + \frac{1}{e}\right)\frac{T}{C_2} $ and $C_4 = 1/C_2$.
\end{proof}
\section{Proof of key auxillary lemmas}\label{sec:appendix}
In this section, we provide  proofs of the key technical lemmas stated in Section \ref{sec:prelim}. We design Lyapunov functionals via appropriate scalings to obtain bounds from the paired inequalities.

\subsection{Proof of \lem{lem:ab-aux}}
Recall the paired inequalities \eqref{ab-aux-sys}:
\begin{equation}\label{ab-aux-sys-p}
     \begin{cases}
         a'(t) \leq b(t), \quad a(0) = a_0,\\[2mm]
         b'(t) \leq C\langle t\rangle ^{-\beta}a(t)+g(t), \quad b(0) = b_0.
     \end{cases}
\end{equation}
We consider the three cases on different choices of $\beta$ seperately.

\medskip\noindent
{\bf (i). $\beta > 2$.} We define a Lyapunov functional
\[
\cL_1(t) := \frac{C}{\beta-1}\langle t\rangle^{-\beta+1}a(t) + b(t), \quad \forall\, t\geq 0.
\]
By \eqref{ab-aux-sys-p} we have 
\begin{align*}
\cL_1'(t) & \le \frac{C}{\beta-1}\langle t\rangle^{-\beta+1}a'(t) -C \langle t\rangle^{-\beta}a(t) + b'(t)\\
& \le \frac{C}{\beta-1}\langle t\rangle^{-\beta+1}b(t) -C \langle t\rangle^{-\beta}a(t) + C \langle t\rangle^{-\beta}a(t) +g(t)\leq \frac{C}{\beta-1}\langle t\rangle^{-\beta+1}\cL_1(t)+g(t).
\end{align*}
Applying Gr\"onwall's inequality, we get
\begin{align*}
\cL_1(t) & \le \left(\cL_1(0)+\int_0^t g(s)\ds\right) \exp\!\left(\frac{C}{\beta-1}\int_0^\infty\langle s\rangle^{-\beta+1}\ds\right)\\
& \le \left(\frac{C}{\beta-1}a_0 + b_0 + \int_0^t g(s)\ds\right)\exp\!\Big(\frac{C}{(\beta-1)(\beta-2)}\Big)=:C_1 + C_2\int_0^t g(s)\ds.
\end{align*}
It follows that for any $t\geq0$
\[
b(t)\leq \cL_1(t) \leq C_1 + C_2\cG_1(t),\quad \text{where}\quad \cG_1(t):=\int_0^t g(s)\ds.\]
Substituting this bound into \eqref{ab-aux-sys-p}$_1$ gives
\[
a(t)\le  a_0 + C_1 t + C_2\int_0^t \cG_1(s)\ds,\quad \forall\, t\ge 0,
\]
which yields the desired bounds in (i).

\medskip\noindent
{\bf (ii). $\beta \in [0,2)$.} 
We define the rescaled quantities
\[
\ta(t) := \bar C\gamma\,e^{-\bar C \,\lan t\ran^{\gamma}} a(t), \quad\tb(t) := \lan t\ran^{1-\gamma} \,e^{-\bar C \,\lan t\ran^{\gamma}}b(t),\quad \forall\, t\geq 0,
\]
where the parameters $\gamma\in(0,1]$ and $\bar C>0$ will be chosen later in \eqref{eq:gammaC}.

We first consider the case when there is no source term, namely $g(t)\equiv0$.
By \eqref{ab-aux-sys-p} we have
\begin{align*}
 \ddt \ta(t) & \leq -  (\bar C\gamma)^2 \lan t\ran^{-(1-\gamma)}\,e^{-\bar C \,\lan t\ran^{\gamma}} a(t) + \bar C\gamma\,e^{-\bar C \,\lan t\ran^{\gamma}} b(t) = - \bar C\gamma \,\lan t\ran^{-(1-\gamma)}\big(\ta(t) -\tb(t)\big); \quad \text{and}\\[2mm]
 \ddt \tb(t) & \leq \big((1-\gamma)\,\lan t\ran^{-\gamma} - \bar C\gamma\big) \,e^{-\bar C \,\lan t\ran^{\gamma}}b(t) + \lan t\ran^{1-\gamma} \,e^{-\bar C \,\lan t\ran^{\gamma}} C\,\langle t\rangle ^{-\beta} a(t)\\
    & \leq (1-\gamma -\bar C\gamma)\,\lan t\ran^{-(1-\gamma)}\,\tb(t) + \frac{C}{\bar C\gamma}\lan t\ran^{1-\gamma-\beta}\ta(t),
\end{align*} 
where in the last inequality, we have used the fact that $\lan t\ran^{-\gamma}\leq1$. Now we choose $\gamma$ and $\bar C$ such that $-(1-\gamma -\bar C\gamma)=\frac{C}{\bar C\gamma}$ and $-(1-\gamma) = 1-\gamma-\beta$, namely:
\begin{equation}\label{eq:gammaC}
\gamma:=1-\frac\beta2,\quad\text{and}\quad\bar C:= \dfrac{1-\gamma + \sqrt{(1-\gamma)^2 + 4C}}{2\gamma}.
\end{equation}
Then we deduce that
\[
 \ddt \tb(t)\leq \frac{C}{\bar C\gamma}\,\lan t\ran^{-(1-\gamma)}\,\big(\ta(t) -\tb(t)\big).
\]
Define an Lyapunov functional 
\[
\cL_2(t) := C\ta(t)+(\bar C\gamma)^2\tb(t), \quad \forall\, t\geq 0.
\]
Then we have the monotone property
\[
\cL_2'(t)\leq - C\bar C\gamma \,\lan t\ran^{-(1-\gamma)}\big(\ta(t) -\tb(t)\big) + C\bar C\gamma \,\lan t\ran^{-(1-\gamma)}\big(\ta(t) -\tb(t)\big) = 0.
\]
It yields that $\cL_2(t)\leq \cL_2(0) = C \bar C \gamma e^{-\bar C}a_0 + (\bar C \gamma)^2e^{-\bar C} b_0$, and consequently
\begin{align*}
a(t) & = (\bar C\gamma)^{-1} e^{\bar C \,\lan t\ran^{\gamma}} \ta(t)\leq (\bar C\gamma)^{-1} e^{\bar C \,\lan t\ran^{\gamma}}C^{-1}\cL_2(0) = \big(a_0 + C^{-1}\bar C \gamma b_0\big)\,e^{\bar C (\lan t\ran^{\gamma}-1)};\quad \text{and}\\
b(t) & = \lan t\ran^{-(1-\gamma)} \,e^{\bar C \,\lan t\ran^{\gamma}} \tb(t) \leq \lan t\ran^{-(1-\gamma)} \,e^{\bar C \,\lan t\ran^{\gamma}} (\bar C \gamma)^{-2}\cL_2(0) = \big(C \bar C^{-1} \gamma^{-1} a_0 + b_0\big)\, \lan t\ran^{-(1-\gamma)} \,e^{\bar C (\lan t\ran^{\gamma}-1)}.
\end{align*}

Next, we take into the consideration of the source term $g(t)$. For the Lyapunov functional $\cL_2$, we obtain
\[
\cL_2'(t)\leq (\bar C\gamma)^2\lan t\ran^{1-\gamma} \,e^{-\bar C \,\lan t\ran^{\gamma}} g(t).
\]
Then we have
\[
\cL_2(t)\leq \cL_2(0) + (\bar C\gamma)^2\int_0^t \lan s\ran^{1-\gamma} \,e^{-\bar C \,\lan s\ran^{\gamma}} g(s) \ds = C \bar C e^{-\bar C}\gamma a_0 + (\bar C \gamma)^2 \big(e^{-\bar C}b_0+\cG_2(t)\big),
\]
where we denote the integrand by $\cG_2(t)$. 
We conclude with our desired bounds:
\begin{align*}
a(t) & \leq \big(e^{-\bar C}a_0 + C^{-1}\bar C \gamma (e^{-\bar C}b_0+\cG_2(t))\big)\,e^{\bar C\,\lan t\ran^{\gamma}};\quad \text{and}\\
b(t) & \leq \big(C \bar C^{-1} e^{-\bar C}\gamma^{-1} a_0 + e^{-\bar C}b_0+\cG_2(t)\big)\, \lan t\ran^{-(1-\gamma)} \,e^{\bar C\,\lan t\ran^{\gamma}},
\end{align*}
where $C_3 := (a_0 + C^{-1}\bar C \gamma b_0)e^{-\bar C}$ and $C_4 := (C \bar C^{-1} \gamma^{-1} a_0 + b_0) e^{-\bar C}$.

\medskip\noindent
{\bf (iii). $\beta = 2$.} This case can be treated similar to the proof of (ii). We define the rescaled quantities
\[
\ta(t):= \lan t\ran^{-\zeta}a(t), \quad \tb(t):= \zeta^{-1}\,\lan t\ran^{-(\zeta-1)}b(t),\quad \forall\, t\geq 0.
\]
By \eqref{ab-aux-sys-p}, we have:
\begin{align*}
    \ddt\ta(t) & \leq - \zeta \lan t\ran^{-\zeta-1} a(t) + \lan t\ran^{-\zeta}b(t) = \zeta\, \lan t\ran^{-1}\big(\tb(t)-\ta(t)\big);\quad \text{and} \\[2mm]
    \ddt \tb(t) & \leq -\zeta^{-1}(\zeta -1)\lan t\ran^{-\zeta}b(t) +  C\zeta^{-1}\lan t\ran^{-(\zeta+1)} a(t) + \zeta^{-1}\,\lan t\ran^{-(\zeta-1)}g(t)\\
    & = C\zeta^{-1} \lan t\ran^{-1}\big(\ta(t) -\tb(t)\big) +  C\zeta^{-1}\lan t\ran^{-(\zeta+1)} a(t) + \zeta^{-1}\,\lan t\ran^{-(\zeta-1)}g(t).
\end{align*} 
where in the last equality, we have used $\zeta( \zeta - 1) = C$ with our choice of $\zeta= (1+\sqrt{1+4C})/2$. 

Define the Lyapunov functional $\cL_3(t):= C\ta(t) + \zeta^2\tb(t)$ so that $\cL_3'(t)\leq \zeta\,\lan t\ran^{-(\zeta-1)}g(t)$, and hence
\[
\cL_3(t)\leq \cL_3(0) + \zeta\int_0^t \lan s\ran^{-(\zeta-1)} \, g(s) \ds = C a_0 + \zeta \big(b_0+\cG_3(t)\big),\quad\text{where}\quad \cG_3(t):=\int_0^t \lan s\ran^{1-\gamma} \, g(s) \ds.
\] 
It then yields
\begin{align*}
 a(t) & = \lan t\ran^{\zeta}\, \ta(t)\leq \big(a_0+C^{-1}\zeta (b_0+\cG_3(t))\big)\,\lan t \ran^{\zeta};\quad \text{and}\\
 b(t) & = \zeta\,\lan t\ran^{\zeta-1} \tb(t)\leq \big(C\zeta^{-1}a_0+b_0+\cG_3(t)\big) \lan t\ran^{\zeta-1},
\end{align*}
finishing the proof of (iii) with $C_5 = a_0+C^{-1}\zeta b_0$ and $C_6 = C\zeta^{-1}a_0+b_0$.
\subsection{Proof of \lem{lem:ab-aux2}}
Recall the paired inequalities \eqref{ab-aux-sys2}
\begin{equation}\label{ab-aux-sys2-p}
     \begin{cases}
         a'(t)\leq b(t), \quad a(0) = a_0,\\[2mm]
         b'(t) \leq C\langle t\rangle ^{-2} (\log \lan t\ran )^{\frac{2\alpha}{1-\alpha}} a(t) + g(t), \quad b(0) = b_0.
     \end{cases}
\end{equation}
When $\alpha=0$, the system reduces to \eqref{ab-aux-sys-p} with $\beta=2$. The result follows from Lemma \ref{lem:ab-aux} (iii). We focus on the proof of the case when $\alpha\in(0,1)$. Due to the additional logarithmic term in \eqref{ab-aux-sys2-p}$_2$, we define a new set of rescaled quantities
 \[
 \ta(t) = \bar{C}\theta e^{-\bar C(\log \langle t\rangle)^{\theta}}a(t), \quad \tb(t) = \langle t\rangle(\log \langle t\rangle)^{-(\theta-1)}\,e^{-\bar C (\log \langle t\rangle)^{\theta}}b(t), \quad \forall\,t\geq 0,
 \]
 where the parameters $\theta>1$ and $\bar C>0$ will be chosen later in \eqref{eq:thetaC}.
 By \eqref{ab-aux-sys2-p}, we have
 \begin{align*}
     \ddt \ta(t) & \leq  \bar{C}\theta e^{-\bar C(\log \langle t\rangle)^{\theta}} b(t) - \bar{C}\theta\, \langle t\rangle^{-1}\,(\log \langle t\rangle)^{\theta-1}\,e^{-\bar C (\log \langle t\rangle)^{\theta}}a(t)\\
     & \leq  \bar{C}\theta\, \langle t\rangle^{-1}\,(\log \langle t\rangle)^{\theta-1}\big(\tb(t) - \ta(t)\big),\\[2mm]
     \ddt \tb(t) & \leq \left((\log \langle t\rangle)^{-(\theta-1)}-(\theta-1) (\log \langle t\rangle)^{-\theta} - \bar C \theta\right)\,e^{-\bar C (\log \langle t\rangle)^{\theta}}b(t)\\
     & \qquad + \langle t\rangle(\log \langle t\rangle)^{-(\theta-1)}\,e^{-\bar C (\log \langle t\rangle)^{\theta}} \left(C\langle t\rangle ^{-2} (\log \lan t\ran )^{\frac{2\alpha}{1-\alpha}} a(t)+g(t)\right)\\
     & \leq \big(1-\bar C\theta \big)\,\lan t\ran^{-1}(\log \lan t\ran )^{\theta-1}\tb(t) + \frac{C}{\bar{C}\theta}\lan t\ran^{-1}(\log \lan t\ran )^{-(\theta-1)+\frac{2\alpha}{1-\alpha}}\,\ta(t)\\
     & \qquad +\langle t\rangle(\log \langle t\rangle)^{-(\theta-1)}\,e^{-\bar C (\log \langle t\rangle)^{\theta}}g(t).
 \end{align*}
 Now we choose $\theta$ and $\bar C$ such that $\theta-1 = -(\theta-1)+\frac{2\alpha}{1-\alpha}$ and $\bar C\theta - 1 = \frac{C}{\bar C\theta}$, namely:
 \begin{equation}\label{eq:thetaC}
  \theta:=\frac{1}{1-\alpha},\quad\text{and}\quad \bar C:=\frac{1 +\sqrt{1 + 4C}}{2\theta}.	
 \end{equation}
Then we deduce that
\[
 \ddt \tb(t)\leq  C (\bar C\theta)^{-1} \lan t\ran^{-1}(\log \lan t\ran )^{\theta-1}\big(\ta(t)-\tb(t)\big) + \langle t\rangle(\log \langle t\rangle)^{-(\theta-1)}\,e^{-\bar C (\log \langle t\rangle)^{\theta}}g(t).
\] 
Define a Lyapunov functional
\[
\cL_4(t):= C\ta(t)+(\bar C\theta)^2\tb(t), \quad \forall\, t\geq 0.
\]
Then adding these above inequalities we get 
\[
 \cL_4'(t) \leq (\bar C\theta)^2\, \langle t\rangle(\log \langle t\rangle)^{-(\theta-1)}\,e^{-\bar C (\log \langle t\rangle)^{\theta}}g(t),
\]
which implies
\[
 \cL_4(t) \leq \cL_4(0) + (\bar C\theta)^2\int_0^t \langle s\rangle(\log \langle s\rangle)^{-(\theta-1)}\,e^{-\bar C (\log \langle s\rangle)^{\theta}}g(s) \ds = C\bar C\theta e^{-\bar C}a_0+\big(\bar C\theta)^2 (e^{-\bar C} b_0 +\cG_4(t)\big),
\]
where the integrand is denoted by $\cG_4$. It then yields
\begin{align*}
 a(t) & = (\bar{C}\theta)^{-1} e^{\bar C(\log \langle t\rangle)^{\theta}}\, \ta(t)\leq \big(e^{-\bar C}a_0+C^{-1}\bar C\theta (e^{-\bar C}b_0+\cG_4(t))\big)\,e^{\bar C(\log \langle t\rangle)^{\theta}};\quad \text{and}\\
 b(t) & = \langle t\rangle^{-1}(\log \langle t\rangle)^{\theta-1}\,e^{\bar C (\log \langle t\rangle)^{\theta}} \tb(t)\\
 & \leq \big(C \bar C^{-1} \theta^{-1}e^{-\bar C}a_0+e^{-\bar C}b_0+\cG_3(t)\big) \langle t\rangle^{-1}(\log \langle t\rangle)^{\theta-1}\,e^{\bar C (\log \langle t\rangle)^{\theta}},
\end{align*}
finishing the proof with $C_7=(a_0+C^{-1}\bar C\theta b_0)e^{-\bar C}$ and $C_8=(C \bar C^{-1} \theta^{-1} a_0 + b_0)e^{-\bar C}$.


\begin{thebibliography}{10}

\bibitem{albi2019vehicular}
Giacomo Albi, Nicola Bellomo, Luisa Fermo, S-Y Ha, Jeongho Kim, Lorenzo Pareschi, David Poyato, and Juan Soler.
\newblock Vehicular traffic, crowds, and swarms: From kinetic theory and multiscale methods to applications and research perspectives.
\newblock {\em Mathematical Models and Methods in Applied Sciences}, 29(10):1901--2005, 2019.

\bibitem{ambrosio2021lectures}
Luigi Ambrosio, Elia Bru{\'e}, Daniele Semola, et~al.
\newblock {\em Lectures on optimal transport}, volume 130.
\newblock Springer, 2021.

\bibitem{bahouri2011fourier}
Hajer Bahouri, Jean-Yves Chemin, and Rapha{\"e}l Danchin.
\newblock {\em {F}ourier analysis and nonlinear partial differential equations}, volume 343.
\newblock Springer, 2011.

\bibitem{ben2005opinion}
Eli Ben-Naim.
\newblock Opinion dynamics: rise and fall of political parties.
\newblock {\em Europhysics Letters}, 69(5):671, 2005.

\bibitem{black2024asymptotic}
McKenzie Black and Changhui Tan.
\newblock Asymptotic behaviors for the compressible {E}uler system with nonlinear velocity alignment.
\newblock {\em Journal of Differential Equations}, 380:198--227, 2024.

\bibitem{black2025hydrodynamic}
McKenzie Black and Changhui Tan.
\newblock Hydrodynamic limit of a kinetic flocking model with nonlinear velocity alignment.
\newblock {\em Kinetic and Related Models}, 18(4):609--632, 2025.

\bibitem{bolley2011stochastic}
Fran{\c{c}}ois Bolley, Jos{\'e}~A Canizo, and Jos{\'e}~A Carrillo.
\newblock Stochastic mean-field limit: non-{L}ipschitz forces and swarming.
\newblock {\em Mathematical Models and Methods in Applied Sciences}, 21(11):2179--2210, 2011.

\bibitem{carrillo2014local}
Jos{\'e}~A Carrillo, Young-Pil Choi, and Maxime Hauray.
\newblock Local well-posedness of the generalized {C}ucker-{S}male model with singular kernels.
\newblock {\em ESAIM: Proceedings and Surveys}, 47:17--35, 2014.

\bibitem{choi2025alignment}
Young-Pil Choi, Michal Fabisiak, and Jan Peszek.
\newblock Alignment with nonlinear velocity couplings: Collision avoidance and micro-to-macro mean-field limits.
\newblock {\em SIAM Journal on Mathematical Analysis}, 57(5):5791--5820, 2025.

\bibitem{choi2019cucker}
Young-Pil Choi and Samir Salem.
\newblock {C}ucker-{S}male flocking particles with multiplicative noises: Stochastic mean-field limit and phase transition.
\newblock {\em Kinetic and Related Models}, 12(3):573--592, 2019.

\bibitem{cucker2007emergent}
Felipe Cucker and Steve Smale.
\newblock Emergent behavior in flocks.
\newblock {\em IEEE Transactions on automatic control}, 52(5):852--862, 2007.

\bibitem{dobrushin1979vlasov}
Roland~L'vovich Dobrushin.
\newblock Vlasov equations.
\newblock {\em Functional Analysis and Its Applications}, 13(2):115--123, 1979.

\bibitem{friesen2020stochastic}
Martin Friesen and Oleksandr Kutoviy.
\newblock Stochastic {C}ucker-{S}male flocking dynamics of jump-type.
\newblock {\em Kinetic and Related Models}, 13(2):211--247, 2020.

\bibitem{golse2016dynamics}
Fran{\c{c}}ois Golse.
\newblock On the dynamics of large particle systems in the mean field limit.
\newblock In {\em Macroscopic and large scale phenomena: coarse graining, mean field limits and ergodicity}, pages 1--144. Springer, 2016.

\bibitem{ha2010emergent}
Seung-Yeal Ha, Taeyoung Ha, and Jong-Ho Kim.
\newblock Emergent behavior of a {C}ucker-{S}male type particle model with nonlinear velocity couplings.
\newblock {\em IEEE Transactions on Automatic Control}, 55(7):1679--1683, 2010.

\bibitem{ha2009simple}
Seung-Yeal Ha and Jian-Guo Liu.
\newblock A simple proof of the cucker-smale flocking dynamics and mean-field limit.
\newblock {\em Communications in Mathematical Sciences}, 7(2):297--325, 2009.

\bibitem{ha2008particle}
Seung-Yeal Ha and Eitan Tadmor.
\newblock From particle to kinetic and hydrodynamic descriptions of flocking.
\newblock {\em Kinetic and Related Models}, 1(3):415--435, 2008.

\bibitem{jabin2014review}
Pierre-Emmanuel Jabin.
\newblock A review of the mean field limits for {V}lasov equations.
\newblock {\em Kinetic and Related models}, 7(4):661--711, 2014.

\bibitem{kim2020complete}
Jong-Ho Kim and Jea-Hyun Park.
\newblock Complete characterization of flocking versus nonflocking of {C}ucker--{S}male model with nonlinear velocity couplings.
\newblock {\em Chaos, Solitons and Fractals}, 134:109714, 2020.

\bibitem{motsch2014heterophilious}
Sebastien Motsch and Eitan Tadmor.
\newblock Heterophilious dynamics enhances consensus.
\newblock {\em SIAM review}, 56(4):577--621, 2014.

\bibitem{mucha2018cucker}
Piotr~B Mucha and Jan Peszek.
\newblock The {C}ucker--{S}male equation: singular communication weight, measure-valued solutions and weak-atomic uniqueness.
\newblock {\em Archive for Rational Mechanics and Analysis}, 227:273--308, 2018.

\bibitem{natalini2021mean}
Roberto Natalini and Thierry Paul.
\newblock On the mean field limit for {C}ucker-{S}male models.
\newblock {\em Discrete and Continuous Dynamical Systems-Series B}, 2021.

\bibitem{nguyen2022propagation}
Vinh Nguyen and Roman Shvydkoy.
\newblock Propagation of chaos for the {C}ucker-{S}male systems under heavy tail communication.
\newblock {\em Communications in Partial Differential Equations}, 47(9):1883--1906, 2022.

\bibitem{shvydkoy2021dynamics}
Roman Shvydkoy.
\newblock {\em Dynamics and analysis of alignment models of collective behavior}.
\newblock Springer, 2021.

\bibitem{shvydkoy2024environmental}
Roman Shvydkoy.
\newblock Environmental averaging.
\newblock {\em EMS Surveys in Mathematical Sciences}, 11(2):277--413, 2024.

\bibitem{sznitman2006topics}
Alain-Sol Sznitman.
\newblock Topics in propagation of chaos.
\newblock In {\em Ecole d'{\'e}t{\'e} de probabilit{\'e}s de Saint-Flour XIX---1989}, pages 165--251. Springer, 2006.

\bibitem{tadmor2023swarming}
Eitan Tadmor.
\newblock Swarming: hydrodynamic alignment with pressure.
\newblock {\em Bulletin of the American Mathematical Society}, 60(3):285--325, 2023.

\bibitem{vicsek2012collective}
Tam{\'a}s Vicsek and Anna Zafeiris.
\newblock Collective motion.
\newblock {\em Physics reports}, 517(3-4):71--140, 2012.

\bibitem{wen2012flocking}
Guanghui Wen, Zhisheng Duan, Zhongkui Li, and Guanrong Chen.
\newblock Flocking of multi-agent dynamical systems with intermittent nonlinear velocity measurements.
\newblock {\em International Journal of Robust and Nonlinear Control}, 22(16):1790--1805, 2012.

\end{thebibliography}
\end{document}